\pdfoutput=1
\documentclass[11pt]{article}
\usepackage{amsfonts}
\usepackage{amsmath}
\usepackage{amssymb}
\usepackage{amsthm}
\usepackage[usenames,dvipsnames]{color}
\usepackage{natbib}{\tiny }
\usepackage[ruled, lined, linesnumbered]{algorithm2e} % Must load after natbib! Else strange errors.
\usepackage{tikz}
\usepackage[utf8]{inputenc}

\usepackage{graphicx}
\usepackage{xfrac}
\usepackage{thmtools, thm-restate}
\usepackage{cancel}
\usepackage{bbm}
\usepackage{bm}
\usepackage[sf,bf,tiny]{titlesec}
\titlelabel{\thetitle. \enspace}

\theoremstyle{definition}
\newtheorem{thm}{Theorem}

\newtheorem{cor}{Corollary}

\def\bbN{\mathbb{N}}

\def\gq{\theta}

\def\bfmath#1{\boldsymbol{#1}}

\def\bfk{{\bfmath{k}}}

\def\ds{\ensuremath{\displaystyle}}

\title{\Large \bf Wright-Fisher diffusion bridges}

\author{
\sc Robert Griffiths\thanks{Department of Statistics, University of Oxford; email: griff@stats.ox.ac.uk}
\and
\sc Paul A.\ Jenkins\thanks{Department of Statistics and Department of Computer Science, University of Warwick; email: P.Jenkins@warwick.ac.uk}
\and 
\sc Dario Span{\`o}\thanks{Department of Statistics, University of Warwick; email: D.Span\`o@warwick.ac.uk }
}

\begin{document}
\maketitle
\noindent
$^*$ Corresponding author.
Version 3.1: TPB 2nd revision
\bigskip

\noindent
{\bf Abstract} 
The trajectory of the frequency of an allele which begins at $x$ at time $0$ and is known to have frequency $z$ at time $T$ can be modelled by the bridge process of the Wright-Fisher diffusion. Bridges when $x=z=0$ are particularly interesting because they model the trajectory of the frequency of an allele which appears at a time, then is lost by random drift or mutation after a time $T$. The coalescent genealogy back in time of a population in a neutral Wright-Fisher diffusion process is well understood. In this paper we obtain a new interpretation of the coalescent genealogy of the population in a bridge from a time $t\in (0,T)$. In a bridge with allele frequencies of 0 at times 0 and $T$ the coalescence structure is that the population coalesces in two directions from $t$ to $0$ and $t$ to $T$ such that there is just one lineage of the allele under consideration at times $0$ and $T$. 
The genealogy in Wright-Fisher diffusion bridges with selection is more complex than in the neutral model, but still with the property of the population branching and coalescing in two directions from time $t\in (0,T)$. The density of the frequency of an allele at time $t$ is expressed in a way that shows coalescence in the two directions. 
A new algorithm for exact simulation of a neutral Wright-Fisher bridge is derived. This follows from knowing the density of the frequency in a bridge and exact simulation from the Wright-Fisher diffusion. The genealogy of the neutral Wright-Fisher bridge is also modelled by branching P\'olya urns, extending a representation in a Wright-Fisher diffusion. This is a new very interesting representation that relates Wright-Fisher bridges to classical urn models in a Bayesian setting.

This paper is dedicated to the memory of Paul Joyce.
%%%%%%%%%%%%%%%%%%%%%%%%%%%
%This paper considers distribution theory and coalescent genealogy in a Wright-Fisher diffusion process, modelling the frequency of an allele, when the frequency is fixed in a bridge at two times $0$ and $T$. Particular attention is paid to beginning and ending at zero frequency. Bridges are invariant under Doob $h$-transforms. We find alternative forms for the bridge density and and derive new identities for coalescent transition functions from these identities.
%The coalescent genealogy of the population at time $t$ is described by two coalescent processes from $t$ to $0$ and from $t$ to $T$. This description of genealogy in bridges is new.
%Genealogy of the neutral Wright-Fisher bridge is also modelled in a new way by branching P\'olya urns, extending a representation in a Wright-Fisher diffusion. 
%%
% A new algorithm for exact simulation of a neutral Wright-Fisher bridge is derived. This follows from exact simulation algorithms for the Wright-Fisher density and an identity for the bridge density from $h$-transformation.
%The genealogy in a Wright-Fisher diffusion with selection is more complex than in a neutral model. It is a branching coalescing random graph. We consider new theory and a coalescent interpretation analogous to the neutral model bridge theory as much as possible.
%
%This paper is dedicated to the memory of Paul Joyce.
\bigskip
%\noindent
%{\bf MSC} 

\noindent
{\bf Keywords:} Wright-Fisher diffusion bridges, Coalescent processes in Wright-Fisher diffusion bridges.

%\noindent
%{\bf Running Head:} 

%%
\section{Introduction}
The trajectory of the frequency of an allele which begins at $x$ at time $0$ and is known to have frequency $z$ at time $T$ can be modelled by the bridge process of the Wright-Fisher diffusion. Henceforth we call this process the Wright-Fisher bridge.
Developing a theory of Wright-Fisher bridges is interesting from a population genetics point of view, particularly when an allele frequency is conditioned to begin at $0$ and end at $0$. This interest is because, informally, in populations where there is no mutation there is a recurrent process where individual genes begin with $0$ frequency and are either lost at $0$ frequency or fix at frequency $1$ at a subsequent time. This also occurs in the infinitely-many-alleles model where new mutants arise at frequency $0$ and are lost at frequency $0$ because of mutation away from any particular type. Wright-Fisher bridges are also interesting from a theoretical point of view  as to how they fit into the general theory of bridges in diffusion processes. The Kingman coalescent and variations, \citep{G1980,K1982,T1984, BEG2000}, look backward in time in Wright-Fisher diffusion processes and are important in describing genealogy in populations.
They are formally dual processes to the diffusion processes. 
\citet{SGE2013} investigated bridges in a Wright-Fisher process in a neutral model and a model with selection. 
 The neutral model genealogy plays a role in obtaining theoretical results in \citet{SGE2013}. The main emphasis in their paper was a clever rejection sampling simulation of the bridge beginning and ending with zero frequency when there is genic selection acting on the allele in the model. In the current paper the genealogy in neutral bridges and bridges where there is genic selection is investigated.
Emphasis in this paper is on a thorough understanding and description of coalescent processes in the bridge. The coalescent interpretation in bridges is new.
The neutral model and the model with selection are both reversible before fixation occurs which implies a bridge from time $0$ to $T$ looks the same probabilistically in reverse time from $T$ to $0$.
%The genealogy of bridges from $0$ to $0$ frequency is interesting in that if the frequency of the allele under consideration at time $t$ is $y$, then coalescence of the genes of frequency $y$ must happen in two directions to a single lineage from $t$ to time 0 and from $t$ to time $T$. This is illustrated in Figure 1 in a neutral model and in Figure 4 in a model with selection.
The genealogy in a Wright-Fisher diffusion with selection is more complex than in a neutral model.
It is a branching coalescing random graph.
We consider theory analogous to the neutral model bridge theory as much as possible, using the transition probabilities of \cite{BEG2000} for the branching coalescing genealogy.

%Genealogy of the neutral Wright-Fisher bridge is also modelled by branching P\'olya urns, extending a representation in a Wright-Fisher diffusion by \citet{GS2010}.
%This is a new very interesting representation that relates Wright-Fisher bridges to classical urn models in a Bayesian setting.

Section \ref{stwo:0} describes known results about Wright-Fisher diffusion processes with and without selection that are necessary to understand bridges. We are particularly concerned with transition density expansions which will allow us to make calculations in a bridge probability density and to explain the coalescent genealogy in the bridge. A spectral expansion of the transition density by \cite{K1964} is well known, however an expansion as a mixture in terms of transition functions in a dual coalescent process is important and less well known.

Section \ref{chap3:0} contains the main results on Wright-Fisher bridges. 
We begin by explaining the use of Doob $h$-transforms in the derivation of the density of an allele frequency in a bridge.
%We begin by explaining the Doob $h$-transform and how the density in a bridge is invariant under $h$-transformation of the unconditioned process. An application of the appropriate $h$-transform in a Wright-Fisher diffusion converts the transition density in a model with no mutation, or mutation away from the allele being considered to the transition density in a model with recurrent mutation. This is important in providing an alternative expression for the bridge density where the allele frequency begins at $0$, since the transformed unconditioned transition density is then proper. The alternative form is also important in exact simulation from the bridge.
In Section \ref{NWFB} neutral Wright-Fisher bridges when there is no mutation are considered in detail with a new interpretation of coalescent genealogy in the bridge. 
The genealogical form of the transition density of the Wright-Fisher diffusion in Section \ref{2.1} is important in showing the coalescent genealogy in Wright-Fisher bridges. An interpretation of Theorem \ref{thm:b} in Section \ref{NWFB} shows that in a bridge with allele frequency beginning at $0$ and ending at $0$ the coalescence structure of the population from time $t \in (0,T)$ is that the population coalesces in two directions from $t$ to $0$ and $t$ to $T$ such that there is just one lineage of the allele under consideration at times $0$ and $T$. This is illustrated in Figure 1. 
The new case of the infinitely-many-alleles model where mutation is away from the allele type is also considered by finding a similar coalescent interpretation where lineages are now lost by coalescence or mutation.
In Section \ref{3.2} there is a subsection on $h$-transformations that condition on non-absorption of the allele frequency and give rise to new interesting identities for the coalescent transition functions and the Jacobi polynomials in the spectral expansions for the unconditional transition density.
 In Section \ref{3.3} we consider the genealogy of the neutral Wright-Fisher bridge modelled by branching P\'olya urns, extending a representation in a Wright-Fisher diffusion by \citet{GS2010}. This is a new very interesting representation that relates Wright-Fisher bridges to classical urn models in a Bayesian setting.
In Section \ref{3.4} Wright-Fisher bridges when there is selection in the model are considered. The genealogy is a branching coalescing graph more complex than the tree in the neutral model, but still with the property of the population branching and coalescing in two directions from time $t\in (0,T)$. This is illustrated in Figure 4. The form of the bridge distribution and coalescent interpretation are new.
 We consider the $h$-transform of the Wright-Fisher transition density with selection conditioned on non-absorption in the future. The transform is related to the first eigenvalue-eigenfunction pair. The Yaglom density is calculated from the first eigenfunction. This is a spheroidal wave function which we show how to calculate with Mathematica using Wolfram Alpha.

In the Appendix a method of exact simulation is developed for neutral Wright-Fisher bridges. 
The algorithm is based on the coalescent genealogy in the bridge and the alternative $h$-transform bridge density.
It is not easy to do exact simulation so this new algorithm is interesting.
\section{The Wright-Fisher diffusion process}\label{stwo:0}
In this section we describe known results about Wright-Fisher diffusion processes with and without selection that are necessary to understand bridges. Genealogical interpretations which give rise to transition density expansions are important in calculating bridge transition functions. 

Let $A$ and $a$ be two types of a gene in a population of individuals. A Wright-Fisher diffusion process modelling the relative frequency of $A$ genes over time has a generator
\begin{equation}
{\cal L}=
\frac{1}{2}x(1-x)\frac{\partial^2}{\partial x^2}
+ \frac{1}{2}\Bigg (x(1-x)\gamma + (\theta_1-\theta x)\Bigg )\frac{\partial}{\partial x}.
\label{pgen:0}
\end{equation}
The population is subject to genic selection, whose sign and strength are described by $\gamma$; and mutations $a \to A$ and $A \to a$ which occur respectively at rates $\theta_1/2$, $\theta_2/2$, with notation $\theta = \theta_1+\theta_2$.
% The transition density of $X(t)$ given $X(0)=x$ will be denoted by $f(x,y;t)$.
\subsection{Neutral Wright-Fisher diffusion}\label{2.1}
Let $f_{\theta_1,\theta_2}(x,y;t)$ be the transition density of the diffusion process with generator (\ref{pgen:0}) when $\gamma=0$ so there is no selection in the model.
There are two forms of a transition density expansion of $f_{\theta_1,\theta_2}(x,y;t)$. The spectral expansion was derived by \citet{K1964} as 
\begin{equation}
f_{\theta_1,\theta_2}(x,y;t) = {\cal B}_{\theta_1,\theta_2}(y)
%
%\nonumber \\
%&&~~~~~~~~\times
\Big \{1 + \sum_{n=1}^\infty e^{-n(n+\theta -1)t/2}
{\widehat{P}}^{(\theta_1,\theta_2)}_n(x)
{\widehat{P}}^{(\theta_1,\theta_2)}_n(y)
\Big \},
\label{fexp:0}
\end{equation}
for $\theta_1,\theta_2 \geq 0$, where 
\[
{\cal B}_{\theta_1,\theta_2}(y) = B(\theta_1,\theta_2)^{-1}y^{\theta_1-1}(1-y)^{\theta_2-1},\>0 < y < 1
\]
is the Beta $(\theta_1,\theta_2)$ density and $\{{\widehat{P}}^{(\theta_1\theta_2)}_n(y)\}_{n=0}^\infty$ are orthonormal polynomials on the Beta density derived from scaling the usual Jacobi polynomials
 $\{{\widetilde{P}}^{(a,b)}_n(r)\}_{n=0}^\infty$ which are well defined for $a,b \geq -1$ and are orthogonal on 
 $(1-r)^a(1+r)^b$, $-1 < r < 1$.  See, for example \citet{I2005}. An expression for these polynomials is
 \[
 \widetilde{P}_n^{(a,b)}(r) = \frac{(a+1)_{(n)}}{n!}\>
  _2F_1\big (-n,a+b+n+1;a+1;(1-r)/2\big ),
 \]
 where $ _2F_1$ is a hypergeometric function.
 %%%%
 %\[
% P_1^{(\theta-1,0)}(r) = \frac{1}{2}\big (r-1 + \theta(r+1)\big ).
% \]
 %%%%%
 In terms of these polynomials the expansion (\ref{fexp:0}) is 
\begin{eqnarray}
&&f_{\theta_1,\theta_2}(x,y;t) = y^{\theta_1-1}(1-y)^{\theta_2-1}
\Bigg \{\sum_{n=0}^\infty e^{-n(n+\theta -1)t/2}c_n(\theta_1,\theta_2)
\nonumber \\
&&~~~~~~~~~~~~~~~~~~~~~~~~~~~~~\times\widetilde{P}^{(\theta_2-1,\theta_1-1)}_n(r)
\widetilde{P}^{(\theta_2-1,\theta_1-1)}_n(s)
\Bigg \},
\label{fexp:1}
\end{eqnarray}
where $r=2x-1$, $s=2y-1$, and
\[
c_n(\theta_1,\theta_2) =
\frac{n!\Gamma(n+\theta_1+\theta_2-1)(2n+\theta_1+\theta_2-1)}
{\Gamma(n+\theta_1)\Gamma(n+\theta_2)}.
\]
The expansion (\ref{fexp:1}) holds for $\theta_1,\theta_2 \geq 0$ taking care with the starting summation index $n$. If $\theta_1,\theta_2 >0$ then the starting index is $n=0$, $\widetilde{P}_0^{(\theta_2-1,\theta_1-1)}\equiv 1$ and the first term is $c_0(\theta_1,\theta_2) =B(\theta_1,\theta_2)^{-1}$. If one of $\theta_1,\theta_2$ is zero and the other non-zero then the summation begins from $n=1$, and if $\theta_1=\theta_2=0$ then the summation begins from $n=2$. In the last case the expansion found by Kimura was the form
\begin{equation}
f_{0,0}(x,y;t) =x(1-x) \sum_{i=1}^\infty e^{-i(i+1)t/2}  (2i+1)i(i+1)R_{i-1}(r)R_{i-1}(s).
\label{fexp:2}
\end{equation}
The polynomials $\{R_i(r)\}_{i=0}^\infty$ are scaled Jacobi polynomials $\{\widetilde{P}_i^{(1,1)}(r)\}_{i=0}^\infty$ orthogonal on $x(1-x),\>0  < x < 1$  such that $R_i(1)=1$. There is an identity between these polynomials and the Jacobi polynomials with index $(-1,-1)$ shown later in the paper as (\ref{Pidentity:0}). If $\theta_1,\theta_2 >0$ then $\{X(t)\}_{t\geq 0}$ is reversible with respect to the stationary Beta $(\theta_1,\theta_2)$ distribution; or if at least one of $\theta_1,\theta_2$ is zero $\{X(t)\}_{t\geq 0}$ is reversible, before being fixed, with respect to the speed measure. 

The Yaglom density
\[
\lim_{t\to \infty}\frac{f_{0,0}(x,y;t)}{\int_0^1f_{0,0}(x,y;t)dy}
\]
is the limit density of the diffusion process conditioned on the diffusion not being fixed at $0$ or $1$ in $(0,t)$. The Yaglom density is straightforward to obtain from (\ref{fexp:1}) or (\ref{fexp:2})
 when either or both of $\theta_1,\theta_2$ are zero  as being proportional to the first orthogonal polynomial times the speed measure density. If $\theta_1=0,\theta_2=\theta$ the Yaglom density is
 \[
 \text{constant}\times y \times y^{-1}(1-y)^{\theta-1} = \theta (1-y)^{\theta-1},\>0 < y < 1.
 \]
 If $\theta_1=\theta_2=0$ the density is 
 \[
 \text{constant}\times y(1-y) \times y^{-1}(1-y)^{-1}= 1,\>0 < y < 1.
 \]
				A second form of the transition density depends on the coalescent genealogy of the population. The Kingman coalescent is a death process which counts the number of lineages back in time in a coalescent tree beginning with a finite or infinite number of leaves. The death rates are ${k\choose 2}$, $k=2,3,\ldots$.
Let $L^\theta(t)$ be the number of non-mutant lineages at time $t$ back in a coalescent tree, beginning with $L^\theta(0)$ leaves, which can be finite or infinity. Mutations occur according to a Poisson process of rate $\theta/2$ on the edges of the coalescent tree, so the death rate of non-mutant lineages is 
${k\choose 2}+k\theta/2$, $k=1,2,\ldots.$. The last lineage is lost by mutation. If there is no mutation $L^0(t)$ counts the number of lineages in the Kingman coalescent.
%\textcolor{red}{Do we want to say anything on duality here?}\\
$\{L^\theta(t)\}_{t\geq 0}$ is a dual process to the Wright-Fisher diffusion which describes the behaviour of the population back in time. The duality is described in \cite{EG1993,GS2013}.
Denote the transition functions of $L^\theta(t)$ when $L^\theta(0)=\infty$ by
\[
P(L^\theta(t) = l\mid L^\theta(0)=\infty) = q_l^\theta(t).
\] 
Notation in this paper will be that for an integer $k\geq 0$ and real number $a$, $a_{[k]}=a(a-1)\cdots (a-k+1)$ and $a_{(k)} = a(a+1)\cdots (a+k-1)$.
An explicit expression \citep{G1980, T1984} is
\[
q_l^\theta(t) = \sum_{k=l}^\infty\rho^\theta_k(t)
(-1)^{k-l}
\frac{(2k+\theta-1)(l+\theta)_{(k-1)}}
{l!(k-l)!},
\]
where $\rho_k^\theta(t) = \exp \{-\frac{1}{2}k(k+\theta-1)t\}$.
The falling  factorial moments of $L^\theta(t)$ are known from \citet{G1980,T1984} as
\begin{equation}
\mathbb{E}\big [L^\theta(t)_{[k]}\big ]
=\sum_{l=k}^\infty\rho_l^\theta(t)(2l+\theta-1){l-1\choose k-1}(\theta+l)_{(k-1)}.
\label{ffmoments:0}
\end{equation}
Let $\{T_l\}_{l=1}^\infty$ be times between events when there are $l$ and $l-1$ non-mutant lineages. The density of $\sum_{k=l}^\infty T_k$ is
\begin{equation*}
\varphi_l^\theta(t) = \frac{l(l+\theta -1)}{2}q^\theta(t),\>t > 0.
%\label{densityid:0}
\end{equation*}
The transition density of $X(t)$ can be expressed as a mixture of Beta densities
\begin{equation}
f_{\theta_1,\theta_2}(x,y;t) = 
\sum_{l=0}^\infty q^\theta_l(t)\sum_{k=0}^l{l\choose k}x^k(1-x)^{l-k}{\cal B}_{k+\theta_1,l-k+\theta_2}(y).
%\nonumber \\
%&&~~~~~~~~~~~~~~~~\times B(k+\theta_1,l-k+\theta_2)^{-1}y^{k + \theta_1-1}(1-y)^{l-k+\theta_2-1}.~~~~~
\label{fexp:3}
\end{equation}
The expansion (\ref{fexp:3}) is a special case of a more general model in \citet{EG1993}. It appeared first implicitly in \cite{G1980} and then written formally in \cite{GL1983}. The range of the summation index $k$ depends on whether $\theta_1,\theta_2$ are zero or positive. If $\theta_1,\theta_2>0$ then $0 \leq k \leq l$; if
$\theta_1=0,\theta_2>0$ then $1 \leq k \leq l$; if $\theta_1>0,\theta_2=0$ then $0 \leq k \leq l-1$; and if $\theta_1=0,\theta_2=0$ then $1 \leq k \leq l-1$.
\citet{GS2010} show the algebraic connection between the two forms of the transition density (\ref{fexp:0}) and (\ref{fexp:3}).
The transition density is improper if either of $\theta_1$ or $\theta_2$ are zero because the population may be fixed for one of the allele types by time $t$ when $X(t)=0$ or $X(t)=1$ with positive probability. Then
$
\int_0^1f_{\theta_1,\theta_2}(x,y;t)dy
$
is the probability that the population is not fixed by time $t$.
This probability has a coalescent interpretation found by integrating in (\ref{fexp:3}).

A genealogical interpretation of the Wright-Fisher population when $\theta_1=\theta_2=0$ is that the coalescent tree backward in time of the individuals beginning at $t$ has $l$ founder edges at the origin. The joint distribution of the $l$ family sizes of the founders at time $t$ is Dirichlet $(1,\ldots,1)$. The types of the families at $t$ are determined by the identity of the $l$ edges at time 0. The probability that an edge is of type $A$ is $x$, and $a$ is $1-x$. The second summation in (\ref{fexp:3}) takes account that there are $k$ edges of type $A$ and $l-k$ of type $a$.
The population is not fixed for one of the allele types if and only if $l\geq 1$ and $l-k \geq 1$.
Grouping the family sizes into $k$ and $l-k$ in the Dirichlet distribution gives the distribution of the frequency of $A$ alleles as Beta $(k,l-k)$. 
%This description is illustrated in Figure 1.
If $\theta >0$ then families at time $t$ can be labelled by those from non-mutant edges at time zero and those from each mutation occurring on the coalescent tree. The joint distribution of the family sizes if there are $l$ non-mutant edges at time zero is
\begin{equation*}
V\bm{U}, (1-V)\{w_{(j)}\}
%\label{families:0}
\end{equation*}
where $\bm{U}$ has a $l$-dimensional Dirichlet$(1,1,\ldots,1)$ distribution; $V$ has a Beta $(l,\theta)$ distribution; $\{w_{(j)}\}$ has a Poisson-Dirichlet $(\theta)$ distribution and the three random variables are independent. This genealogical interpretation occurs in \citet{G1980,EG1993}. If there are two types and $\theta_1,\theta_2 >0$ then choosing the type of the non-mutant edges as $A$ with probability $x$ and mutant types to be of type $A$ with probability $\theta_1/\theta$ gives the density (\ref{fexp:3}). If $\theta_1 =0,\theta_2 >0$ the non-mutant edges type is chosen in the same way and all mutations are of type $a$ instead.
\subsection{Wright-Fisher diffusion with selection}
If $\gamma \ne 0$ selection is present in the model. When $\theta_1,\theta_2>0$ the stationary density of $\{X(t)\}_{t\geq 0}$ is a weighted Beta density
\begin{equation*}
c(\bm{\theta})^{-1}e^{\gamma x}B(\theta_1,\theta_2)^{-1}x^{\theta_1-1}(1-x)^{\theta_2-1},\>0<x<1,
%\label{stationary:1}
\end{equation*}
where $c(\bm{\theta}) = \mathbb{E}\big [e^{\gamma X}\big ]$, with expectation in a Beta $(\theta_1,\theta_2)$ distribution.
$\{X(t)\}_{t\geq 0}$ is reversible with respect to the stationary distribution if $\theta_1,\theta_2 >0$ and reversible with respect to the speed measure if one or both of $\theta_1,\theta_2$ are zero. The eigenvalues and eigenfunctions $\{\lambda_j,w_j\}$ satisfy ${\cal L}w_j = \lambda_jw_j$, or
\begin{equation} 
\frac{1}{2}x(1-x)w_j^{\prime\prime}
+ \frac{1}{2}\Bigg (\gamma x(1-x) + (\theta_1-\theta x)\Bigg )w_j^\prime + \lambda_jw_j=0.
\label{eigenfunction:0}
\end{equation}
The eigenfunctions are orthogonal on the weighted Beta distribution or speed measure, but are not polynomials. If $\theta_1,\theta_2>0$ then $\lambda_0=1$ and $w_0 = 1$. There is no explicit expression for the eigenvalues.   \citet{SS2012} expand the eigenfunctions in terms of Jacobi polynomials and then use numerical methods to find a computational solution of the transition density. The Yaglom limit density when either or both of $\theta_1,\theta_2$ are zero is still proportional to the speed measure times the first eigenfunction, which does not have a simple form. If $\theta_1=\theta_2=0$, $w$ can be expressed in terms of a spheroidal wave function. This was known by \cite{K1955}. We consider the expression in more detail later in the paper.
 If $\theta >0$ then the eigenfunctions cannot be expressed as spheroidal wave functions multiplied by a function not depending on the eigenfunction.

A genealogical form of the transition density of the diffusion process with generator (\ref{pgen:0}), derived in \cite{BEG2000}, is
\begin{equation*}
f_{\theta_1,\theta_2,\gamma}(x,y,t) = \sum_{\bm{\alpha}\in {\mathbb{Z}_+^2}}
b^{(\theta_1,\theta_2)}_{\bm{\alpha}}(t,x)\pi[\bm{\alpha}+\bm{\theta}](y),
%\label{seldensity:0}
\end{equation*}
where 
\[
\pi[\bm{\theta}](y) = c(\bm{\theta})^{-1}B(\theta_1,\theta_2)^{-1}e^{\gamma y}y^{\theta_1-1}(1-y)^{\theta_2-1}.
\]
%and $c(\bm{\theta})$ is a constant chosen so $\pi[\bm{\theta}](y)$ integrates to 1 if $\theta_1>0,\theta_2>0$.
The expansion is analogous to (\ref{fexp:0}) in the neutral model. If $\theta_1>0,\theta_2>0$ then
$\pi[\bm{\theta}]$ is the stationary density of $\{X(t)\}_{t\geq 0}$.
$b^{(\theta_1,\theta_2)}_{\bm{\alpha}}(t,x)$ are transition functions in a two-dimensional birth and death process $\bm{\alpha}(t)$ beginning with an entrance boundary at infinity where the relative frequencies of the types  are $x,1-x$. If $\theta_1>0,\theta_2>0$ summation is over $\alpha_1\geq 0,\alpha_2 \geq 0$; if $\theta_1=0,\theta_2 > 0$ summation is over $\alpha_1\geq 1,\alpha_2 \geq 0$; if $\theta_1>0,\theta_2=0$ summation is over $\alpha_1\geq 0,\alpha_2 \geq 1$; and if $\theta_1=0,\theta_2=0$ summation is over $\alpha_1\geq 1,\alpha_2 \geq 1$. $\{\bm{\alpha}(t)\}_{t\geq 0}$ is a dual process to the diffusion process describing coalescent genealogy back in time.
The transition rates of $\{\bm{\alpha}(t)\}_{t\geq 0}$ when $\gamma >0$, $\theta \geq 0$, are
\begin{eqnarray*}
q(\bm{\alpha},\bm{\alpha}-\bm{e}_i) &=& \frac{1}{2}\alpha_i(|\bm{\alpha}|+\theta -1)
\frac{c(\bm{\theta}+\bm{\alpha}-\bm{e}_i)}{c(\bm{\theta}+\bm{\alpha})},\>i=1,2
\nonumber \\
q(\bm{\alpha},\bm{\alpha}+\bm{e}_2) &=& \frac{1}{2}\gamma (\theta_2 + \alpha_2)\frac{\bm{|\alpha|}}{|\bm{\alpha}|+\theta}
\frac{c(\bm{\theta} + \bm{\alpha}+\bm{e}_2)}{c(\bm{\theta} + \bm{\alpha})}
\nonumber \\
q(\bm{\alpha},\bm{\alpha}) &=& -\frac{1}{2}\left (\alpha_2\gamma
+|\bm{\alpha}|(|\bm{\alpha}|-1)\right ).
%\label{brates:1}
\end{eqnarray*}
$b^{(\theta_1,\theta_2)}_{\bm{\alpha}}(t,x)$ satisfies the usual forward equations
\begin{eqnarray*}
&&\frac{\partial}{\partial t}b^{(\theta_1,\theta_2)}_{\bm{\alpha}}(t,x)=-q(\bm{\alpha},\bm{\alpha})b^{(\theta_1,\theta_2)}_{\bm{\alpha}}(t,x) 
\nonumber \\
&&~~+ \sum_{i=1}^2\Bigg (
q(\bm{\alpha}+\bm{e}_i,\bm{\alpha})
b^{(\theta_1,\theta_2)}_{\bm{\alpha}+\bm{e}_i}(t,x)
+
q(\bm{\alpha}-\bm{e}_i,\bm{\alpha})
b^{(\theta_1,\theta_2)}_{\bm{\alpha}-\bm{e}_i}(t,x)
\Bigg ),
%\label{qforward:0}
\end{eqnarray*}
and the backwards equation
\begin{equation*}
\frac{\partial}{\partial t}b^{(\theta_1,\theta_2)}_{\bm{\alpha}}(t,x)={\cal L}b^{(\theta_1,\theta_2)}_{\bm{\alpha}}(t,x),
%\label{qbackward:0}
\end{equation*}
where ${\cal L}$ is the generator (\ref{pgen:0}).
The boundary condition is
\begin{equation*}
b^{(\theta_1,\theta_2)}_{\bm{\alpha}}(0,x) = {|\bm{\alpha}|\choose \bm{\alpha}}x^{\alpha_1}(1-x)^{\alpha_2},
%\label{qboundary:0}
\end{equation*}
the same as in a neutral model.
\cite{KN1997,NK1997} introduced an ancestral selection graph, which is a branching coalescing graph describing the coalescent genealogy in a Wright-Fisher model with selection. In their graph the types of ancestral lineages are not fully determined until the type of the ultimate ancestor is reached. The ancestral tree is then read off from the graph. The process here is related, but different, in that the lineages are all typed in the graph from time $t$ back to time $0$. A helpful detailed description of the genealogy in a pre-limit Moran model and in our process is in \cite{EG2009}.

The frequencies $y_1=y, y_2=1-y$ can be decomposed into independent ${\cal PD}(\theta_1+\alpha_1)$ and 
${\cal PD}(\theta_2+\alpha_2)$ processes $\{w_{1i}\}$, $\{w_{2i}\}$ which are then weighted as $V\{w_{1i}\},(1-V)\{w_{2i}\}$ where $V$ is independent of the two processes and has a density $\pi[\bm{\theta} + \bm{\alpha}](v)$. This representation is then at the level of individual family sizes either beginning from the origin, or new mutations.

\section{Wright-Fisher diffusion bridges}\label{chap3:0}
Let $\{X^{x,z,[0,T]}(t)\}_{0 \leq t \leq T}$ be a Wright-Fisher diffusion process $\{X(t)\}_{t\geq 0}$ conditioned on $X(0)=x$ and $X(T)=z$. Notation is similar to that used in \citet{SGE2013}. 
 We dispense with subscripts $\theta_1,\theta_2, \gamma$ on $f$ in the explanation of $h$-transforms that follows to ease notation.
The transition density of $X(t)^{x,z,[0,T]}$ at $v$ given $X(s)^{x,z,[0,T]}=u$ is clearly
\begin{equation}
f_{x,z,[0,T]}(u,v; s,t) = \frac{f(u,v;t-s)f(v,z;T-t)}{f(u,z;T-s)},
\label{bridgetr;0}
\end{equation}
where $f(x,y;t)$ is the transition density in the diffusion process with generator (\ref{pgen:0}).
There is a very elegant theory of Doob $h$-transforms and the invariance of a bridge distribution under a transform. The reader is referred to \citet{RW2000} for an introduction to $h$-transforms. We give a very brief introduction and explain the relevance to bridges. In a $h$-transform the transition density is mapped to 
\begin{equation}f(x,y;t)\frac{h(y,t)}{h(x,0)}.
\label{ttt:0}
\end{equation}
The function $h \geq 0$ is chosen so that (\ref{ttt:0}) is the transition density of a stochastic process.  
Typical choices of $h$ provide an  interpretation that the transformed process is conditioned on not being absorbed at boundary points in the future. In a Wright-Fisher diffusion without mutation $0$ and $1$ are absorbing boundary points.
In a model where absorption is certain let $\tau$ be the time when this occurs. Let $h(x)$ be the first eigenfunction of the generator satisfying ${\cal L}h(x) = -\lambda_1h(x)$, $h(x) \geq 0$, where $\lambda_1 >0$ is the largest eigenvalue and set
$h(x,t) = h(x)e^{-\lambda_1t}$. Then
\[
\lim_{\tau^*\to \infty}\frac{P_x(\tau >\tau^*)}{h(x,\tau^*)} = \text{constant}.
\]
The transition functions (\ref{ttt:0})
can be thought of as asymptotic transition functions of a process conditioned not to be lost or fixed. Let $\tau^*> t$, then the transition functions of a process conditioned on $\tau > \tau^*$ are 
\[
f(x,y;t)\frac{P_y(\tau > \tau^*-t)}{P_x(\tau > \tau^*)}
\]
which converges to (\ref{ttt:0}) as $\tau^*\to \infty$.
The generator of a process with infintesimal variance $\sigma^2(x)$ and mean $\mu(x)$ which is $h$-transformed by a $h$ not depending on time is
\begin{equation*}
\frac{1}{2}\sigma^2(x)\frac{\partial^2}{\partial x^2} 
+ \Big (\sigma^2(x)\frac{d}{dx}\log h(x) + \mu(x)\Big )\frac{\partial}{\partial x}
%\label{hhh:1}
\end{equation*}
and the stationary distribution is
\[
\frac{h(y)}{\int_0^1h(x)dx}.
\]
A more formal description of conditioning in general diffusion processes is the first theorem in \citet{P1985}.

The probability distribution of bridges is invariant
under a $h$-transform.
This is important for us in finding alternative forms for the density of the frequency of an allele in a bridge, particularly when the frequency is $0$ at times $0$ and $T$.
As in \citet{SGE2013} consider the joint density of the $h$-transformed bridge process at times $0 < t_1 < \cdots < t_n < T$. The density is 
\begin{eqnarray}
%\begin{split}
&&\frac{f(x,v_1;t_1)\frac{h(v_1)}{h(x)}f(v_1,v_2;t_2-t_1)\frac{h(v_2)}{h(v_1)}\cdots f(v_n,y;T-t_n)\frac{h(y)}{h(v_n)}}{\frac{h(y)}{h(x)}f(x,y;T)} \nonumber \\
&& \quad = \frac{f(x,v_1;t_1)f(v_1,v_2;t_2-t_1)\cdots f(v_n,y; T-t_n)}{f(x,y;T)}.
%\end{split}
%\]
\label{hBridge:0}
\end{eqnarray}
This shows the invariance of the bridge process under $h$-transforms not depending on time.
In a similar way the bridge density is invariant under a $h$-transform depending on time $h(x,t)=h(x)e^{\lambda t}$, where $\lambda$ is a constant. 
Reversibility of $\{X(t)\}_{t\geq 0}$ implies that the probabilistic behaviour of the bridge process should look the same from $0$ to $T$ as it does backwards in time from $T$ to $0$. Let $\pi(x)$ be the density of the reversing measure so that $\pi(x)f(x,y;t) = \pi(y)f(y,x;t)$. Now consider the $h$-transform with $h(x)=\pi(x)^{-1}$. The transform reverses the path with  (\ref{hBridge:0}) evaluating to
\begin{equation*}
\frac{f(y,v_n;T-t_n)f(v_n,v_{n-1};t_{n-1})\cdots f(v_1,x,t_1)}
{f(y,x;T)},
%\label{hBridge:1}
\end{equation*}
which is the path density back in time from $T$ to $0$ of a bridge from $y$ to $x$.
The next theorem is general when a process is reversible. It appears for a general process in \citet{FPY1992} and in \citet{SGE2013} for the Wright-Fisher diffusion.
\begin{thm}\label{thm:a}
The density of 
$X^{x,z,[0,T]}(t)$
 and the reverse form is
\begin{eqnarray}
f_{x,z,[0,T]}(y;t) &=&\frac{f(x,y;t)f(y,z;T-t)}{f(x,z;T)}
\nonumber \\
&=& \frac{f(z,y;T-t)f(y,x;t)}{f(z,x;T)}
\nonumber \\
&=& f_{z,x,[0,T]}(y;T-t).
\label{bridge:rev:0}
\end{eqnarray}
\end{thm}
\subsection{Neutral Wright-Fisher bridges}\label{NWFB}
In this section the coalescent genealogy of bridges is considered.
 We are particularly interested  in $\{X^{0,0,[0,T]}(t)\}_{0\leq t \leq T}$ bridges. Bridges where there is no mutation are studied, as well as bridges where mutation is away from the allele of interest at a rate $\theta$, as in the infinitely-many-alleles model. An approach is to let $x,z\to 0$ in the general bridge.
Often beginning a diffusion process at $0$ is analytically difficult, but using reversibility and a coalescent genealogy approach the treatment is relatively straightforward.
The transition density as $x \to 0$ is
\begin{eqnarray}
f_{0,0}(x,y;t) &=& 
\sum_{l=2}^\infty q^0_l(t){l\choose 1}x(1-x)^{l-1}
B(1,l-1)^{-1}(1-y)^{l-2} + o(x^2)
\nonumber \\
&=&
x\sum_{l=2}^\infty l(l-1)q^0_l(t)(1-y)^{l-2} + o(x^2).
\label{one:0}
\end{eqnarray}
The term of order $x$ in the first line of (\ref{one:0}) is the joint density of $X(t)$ and the probability that there is exactly one coalescent lineage from time $t$ back to time zero that is of type $A$, and at least one other lineage of type $a$. Denote
\begin{equation}
g(y;t) = \sum_{l=2}^\infty l(l-1)q_l^0(t)(1-y)^{l-2}.
\label{gseries:0}
\end{equation}
An interpretation is that $g(y;t)$ is the probability that there are at least two founder lineages with exactly one lineage of the same type as the allele with frequency $y$. This is because given $l$ founders at time $t$ back the family size at $t$ of a given lineage has density $(l-1)(1-y)^{l-2},\> 0 < y < 1$ and there are $l$ choices for such a lineage.

Note the identity, from (\ref{fexp:3}), that
$f_{0,0}(y,0;t) = y(1-y)g(y;t)$. This is a consequence of reversibility with respect to the speed measure density $y^{-1}(1-y)^{-1}$ because
\begin{eqnarray*}
				\lim_{z\to 0} f_{0,0}(y,z;t) &=& \lim_{z\to 0} \frac{z^{-1}(1-z)^{-1}}{y^{-1}(1-y)^{-1}}f_{0,0}(z,y;t)
\nonumber \\
&=& y(1-y)g(y;t).
\end{eqnarray*}
Now dividing the numerator and denominator of the second line in (\ref{bridge:rev:0}) by $x$,
\begin{equation}
\lim_{x\to 0}f_{x,z,[0,T]}(y;t) = \frac{g(y;t)f(y,z;T-t)}{g(z;T)}.
\label{inter:0}
\end{equation}
Set $z=0$ in (\ref{inter:0}) to obtain that
\begin{equation*}
f_{0,0,[0,T]}(y;t) = \frac{y(1-y)g(y;t)g(y;T-t)}{h(T)},
\end{equation*}
where
\[
h(T) := g(0;T) = \sum_{l=2}^\infty l(l-1)q_l^0(T).
\]
Substituting in (\ref{ffmoments:0}) with $\theta=0$ and $k=2$,
\begin{equation}
h(T)=g(0;T)=
\sum_{l=2}^\infty l(l-1)q_l^0(T) = \sum_{l=2}^\infty e^{-\frac{1}{2}l(l-1)T}
(2l-1)l(l-1).
\label{ffact:0}
\end{equation}
Collecting results and substituting in (\ref{bridge:rev:0}) gives the next theorem.
\begin{thm}\citep{SGE2013}\label{thm:b}
The density of $X^{0,0,[0,T]}(t)$, $0\leq t \leq T$ in a neutral Wright-Fisher diffusion bridge with $\theta_1=\theta_2=0$ is
\begin{equation}
\frac{y(1-y)g(y;t)g(y;T-t)}{\sum_{l=2}^\infty (2l-1)l(l-1)e^{-l(l-1)T/2}}.
\label{thm1:0}
\end{equation}
\end{thm}
There is a clear reversibility in the density (\ref{thm1:0}) from both ends of the bridge.

A diagram of the genealogy of a (0,0) bridge is illustrated in Figure 1. 
\bigskip

\vbox{
\begin{center}
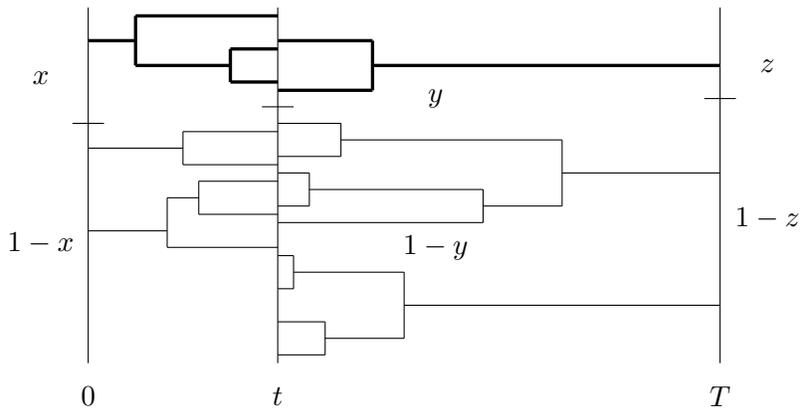
Figure 1: Genealogy in a neutral (0,0) bridge
\bigskip

\begin{tikzpicture}[xscale = 0.21,yscale=0.22]
\draw (0,0) -- (0,21.5);
\draw (12,0) -- (12,15.5);
\draw (12,15.5) -- (12,21.5);
\draw (40,0) -- (40,21.5);
%%%Left side
\draw [line width=1.25pt] (0,19.5) -- (3,19.5);
\draw [line width=1.25pt] (3,21) -- (12,21);
\draw [line width=1.25pt] (3,18) -- (9,18);
\draw [line width=1.25pt] (9,19) -- (12,19);
\draw [line width=1.25pt] (9,17) -- (12,17);
\draw [line width=1.25pt] (3,21) -- (3,18);
\draw [line width=1.25pt] (9,19) -- (9,17);
\draw (0,13) -- (6,13);
\draw (6,14) -- (12,14);
\draw (6,12) -- (12,12);
\draw (6,14) -- (6,12);
\draw (0,8) -- (5,8);
\draw (5,7) -- (12,7);
\draw (5,10) -- (7,10);
\draw (7,11) -- (12,11);
\draw (7,9) -- (12,9);
\draw (5,10) -- (5,7);
\draw (7,11) -- (7,9);
% Right side
\draw [line width=1.25pt] (40,18) -- (18,18);
\draw [line width=1.25pt] (18,19.5) -- (12,19.5);
\draw [line width=1.25pt] (18,16.5) -- (12,16.5);
\draw [line width=1.25pt] (18,19.5) -- (18,16.5);
\draw (40,11.5) -- (30,11.5);
\draw (30,13.5) -- (16,13.5);
\draw (16,14.5) -- (12,14.5);
\draw (16,12.5) -- (12,12.5);
\draw (16,14.5)  -- (16,12.5);
\draw (30,13.5) -- (30,9.5);
\draw (25,8.5) -- (12,8.5);
\draw (25,8.5) -- (25,10.5);
\draw (14,11.5) -- (14,9.5);
\draw (30,9.5) -- (25,9.5);
\draw (25,10.5) -- (14,10.5);
\draw (14,9.5) -- (12,9.5);
\draw (14,11.5) -- (12,11.5);
\draw (40,3.5) -- (20,3.5);
\draw (20,5.5) -- (13,5.5);
\draw (13,4.5) -- (12,4.5);
\draw (13,6.5) -- (12,6.5);
\draw (13,6.5) -- (13,4.5);
\draw (20,1.5) -- (15,1.5);
\draw (15,0.5) -- (12,0.5);
\draw (15,2.5) -- (12,2.5);
\draw (15,2.5) -- (15,0.5);
\draw (20,5.5) -- (20,1.5);
\draw (-1,14.5) -- (1,14.5);
\draw -- (-3,17.25) node {$x$};
\draw -- (-3,7.25) node {$1-x$};
\draw (39,16) -- (41,16);
\draw -- (43,18) node {$z$};
\draw -- (43,8.75) node {$1-z$};
\draw (11,15.5) -- (13,15.5); % Horizontal line for y %
\draw -- (22,16) node {$y$};
\draw -- (22,7) node {$1-y$};
\draw -- (0,-2) node {$0$};
\draw -- (12,-2) node {$t$};
\draw -- (40,-2) node {$T$};
\end{tikzpicture}
\end{center}

There must be exactly 1 lineage from $0$ to $t$ and from $T$ to $t$ in reverse time with leaves in the frequency $y$, conditional on $y$ being reached  as $x,z \to 0$.
}
The description of the coalescent genealogy in the bridge is new in this paper and important.
Informally the bridge density (\ref{thm1:0}) times $dy$ is proportional to the proportion of sample paths in $(y,y+dy)$ at a given time where the allele of frequency $y$ has a single ancestor in both directions at times $t$ and $T-t$, divided by the proportion of sample paths in $(y,y+dy)$, which is the speed measure $y^{-1}(1-y)^{-1}$.

% At a given time $t \in [0,T]$ coalescence of lineages happens in both time directions towards $0$ and $T$. $y(1-y)g(y;t)$ and $y(1-y)g(y;T-t)$ are the probabilities, in a limit as $x,z \to 0$, that there is exactly one coalescent lineage from frequency $y$ hitting the allele frequencies of that type at $0$ and $T$. The density (\ref{thm1:0}) is proportional to the speed measure $y^{-1}(1-y)^{-1}$ times these probabilities.

It is of interest to consider the coalescent genealogy at time points $s < t <T$ and calculate the number of lineages $b$, at $T$, which are non-mutant back to $t$, but not as far back as $s$; the number of non-mutant lineages $a$ back to time $s$; and  the number of distinct non-mutant lineages $c$ back to time $s$. To begin condition on the number of non-mutant lineages from $s$ to $t$ $L^\theta(t-s) = l$.
Then, from the representation
%(\ref{families:0})
at time $t$, considering the roots of the lineages from $T$ back to $t$, it must be that $a$ fall on frequencies from $V\bm{U}$ and $b$ fall on the frequencies $(1-V)\{w_{(j)}\}$. Therefore the conditional probability is
\begin{equation*}
{a+b\choose a}\mathbb{E}\big [V^a(1-V)^b\big ] = {a+b\choose a}
\frac{l_{(a)}\theta_{(b)}}{(l+\theta)_{(a+b)}}.
%\label{vab:0}
\end{equation*}
The number of distinct non-mutant lineages $c$ back to $s$ is the number of distinct frequencies $V\bm{U}$ which are hit by the $a$ root lineages. The required probability is therefore
\begin{equation}
\sum_{a\geq c}q_{a+b}(T-t){a+b\choose a}
\mathbb{E}\big [V^a(1-V)^b\big ]
{l\choose c}\sum_{\bm{a}; a_j \geq 1}
\frac{a!}{a_1!\cdots a_c!}
\mathbb{E}\big [U_1^{a_1}\cdots U_c^{a_c}\big ].
\label{ceq:0}
\end{equation}
Now
\[
\mathbb{E}\big [U_1^{a_1}\cdots U_c^{a_c}\big ]
=\frac{a_1!\cdots a_c!}{l_{(a)}},
\]
so the sum (\ref{ceq:0}) is equal to
\begin{eqnarray}
\sum_{a \geq c}
q_{a+b}(T-t){a+b\choose a}
\frac{l_{(a)}\theta_{(b)}}{(l+\theta)_{(a+b)}}\cdot\frac{a!}{l_{(a)}}{a-1\choose c-1}.
\label{ceq:1}
\end{eqnarray}
The marginal distribution of $c$ from (\ref{ceq:1}) is, setting $m=a+b$,
\begin{equation}
{l\choose c}\sum_{m\geq c}q_m(T-t)
\sum_a\frac{m!}{(m-a)!}\frac{\theta_{(m-a)}}{(l+\theta)_{(m)}}{a-1\choose c-1}.
\label{ceq:2}
\end{equation}
The inner sum of (\ref{ceq:2}) is equal to 
\begin{equation}
\frac{m(1+\theta)_{(m-1)}}{(l+\theta)_{(m)}}\sum_a
{m-1\choose a-1}{a-1\choose c-1}
\frac{1_{(a-1)}\theta_{(m-1-(a-1))}}
{(1+\theta)_{(m-1)}}
\label{ceq:3a}
\end{equation}
which is the coefficient of $u^{c-1}$ in 
\begin{equation}
\frac{m(1+\theta)_{(m-1)}}{(l+\theta)_{(m)}}\sum_{a\geq 1}
{m-1\choose a-1}(1+u)^{a-1}\mathbb{E}\big [\xi^{a-1}(1-\xi)^{m-1-(a-1)}\big ]
\label{ceq:3}
\end{equation}
where $\xi$ has a Beta $(1,\theta)$ distribution. Now (\ref{ceq:3}) evaluates to 
\[
\frac{m(1+\theta)_{(m-1)}}{(l+\theta)_{(m)}}
\mathbb{E}\big [(1+u\xi)^{m-1}\big ]
\]
showing that (\ref{ceq:3}) is equal to
\[
\frac{m(1+\theta)_{(m-1)}}{(l+\theta)_{(m)}}
{m-1\choose c-1}\frac{(c-1)!}{(1+\theta)_{(c-1)}}.
\]
Finally
\begin{eqnarray}
&&P(L^\theta(T-s)=c\mid L^\theta(t-s) = l)
\nonumber \\
&=& {l\choose c}\sum_{m\geq c}q_m^\theta(T-t)\frac{m!}{(m-c)!}
\frac{\Gamma(\theta+m)}{\Gamma(\theta+c)(l+\theta)_{(m)}}.
\label{ceq:4}
\end{eqnarray}
From (\ref{ceq:4}) an identity also follows that
\begin{equation*}
q_c^\theta(T-s) = \sum_lq_l^\theta(t-s)
{l\choose c}\sum_{m\geq c}q_m^\theta(T-t)\frac{m!}{(m-c)!}
\frac{\Gamma(\theta+m)}{\Gamma(\theta+c)(l+\theta)_{(m)}}.
\end{equation*}

We now consider a bridge from 0 to 0 frequency when there is mutation away from a gene of type $A$. The importance of this bridge is that it models a new mutant in the infinitely-many-alleles model which arises at time 0, then is lost at time $T$ with mutation away from the new mutant to other types at rate $\theta/2$. The generator for the frequency of a given allele in a neutral model is
\[
{\cal L}=\frac{1}{2}x(1-x)\frac{\partial^2}{\partial x^2} 
- \frac{1}{2}\theta x\frac{\partial}{\partial x}.
\]
The speed measure density is $y^{-1}(1-y)^{\theta-1}$ and the process is reversible with respect to this measure before absorption. Theorem \ref{thm:a} holds for this model. Again care is needed for a 0,0 bridge.
%Denote the transition density as $f^\theta(x,y;t)$.
As $x \to 0$ 
\begin{equation*}
f_{0,\theta}(x,y;t) = x\sum_{l=1}^\infty q_l^\theta(t)l(l+\theta -1)(1-y)^{l+\theta-2} + o(x^2)
%\label{one:1}
\end{equation*}
and
\begin{equation*}
f_{0,\theta}(y,0;t) = \sum_{l=1}^\infty q_l^\theta(t)l(l+\theta-1)y(1-y)^{l-2}.
%\label{one:1a}
\end{equation*}
Denote
\begin{equation}
g_\theta(y;t) = \sum_{l=1}^\infty l(l+\theta -1)(1-y)^{l+\theta-2}q_l^\theta(t).
\label{gtheta:0}
\end{equation}

A similar calculation to that in Theorem \ref{thm:b}, substituting in the second line of (\ref{bridge:rev:0}) gives the bridge density at $t$.
\begin{thm}\label{thm:c}
The density of $X^{0,0,[0,T]}(t)$, $0\leq t \leq T$ in a neutral Wright-Fisher diffusion bridge with mutation away from $A$ at rate $\theta/2$ is
\begin{equation}
\frac{
y(1-y)^{-\theta+1}g_\theta(y;t)g_\theta(y;T-t)
}
{
\sum_{l=1}^\infty
e^{-l(l+\theta-1)T/2}(2l+\theta-1)l(l+\theta-1)
% \big ((l-1)(l+\theta) + 1\big ) Changed, spotted by Paul - Bob
}.
\label{thm:d}
\end{equation}
\end{thm}
The distribution of the frequency in the bridge in Theorem \ref{thm:c} is the same from both ends as it is when there is no mutation. There is a similar diagram and interpretation to that in Theorem \ref{thm:b}, Figure 1 for the genealogy, except that only non-mutant lines from the frequency $y$ are considered and there can only be single non-mutant founder lineages from each of the frequencies at $0$ and $T$ which have families in the $y$ frequency. Theorem \ref{thm:b} and the coalescent interpretation are new.

The mean frequency in a bridge is found by multiplying (\ref{thm:d}) by $y$ and integrating over $(0,1)$.

\begin{cor} 
\begin{equation*}
\mathbb{E}\big [X^{0,0,[0,T]}(t)\big ]
=
\frac{
\sum_{l,m=1}^\infty 
\frac{l(l+\theta-1)m(m+\theta-1)}{(l+m+\theta-2)(l+m+\theta - 3)}
q_l^{\theta}(t)q_m^{\theta}(T-t)
}
{
\sum_{l=1}^\infty
e^{-l(l+\theta-1)T/2}(2l+\theta-1)\big ((l-1)(l+\theta) + 1\big )
}.
%\label{cor:b}
\end{equation*}
\end{cor}
An algorithm to simulate exactly from the density in a neutral Wright-Fisher bridge is given in the Appendix.

\subsection{$h$-transform identities for bridges}\label{3.2}
There are useful and interesting identities for bridge densities that follow from $h$-transforms. The transforms relate neutral Wright-Fisher transition functions when one or both of $\theta_1,\theta_2$ are zero to Wright-Fisher transition functions with recurrent mutation when $\theta_1,\theta_2>0$. The bridge distribution is invariant under the $h$-transform. Bridges beginning at 0 are easier to deal with after the transformation when beginning at zero, because the Wright-Fisher transition function is then well defined. 

Suppose $\theta_1=0,\theta_2=0$, then an identity is
\begin{equation}
f_{2,2}(x,y;t)=e^tf_{0,0}(x,y;t)y(1-y)/x(1-x).
\label{h2:0}
\end{equation}
 To show the identity (\ref{h2:0}) consider the following calculation with the backward generator.
\begin{eqnarray*}
\frac{\partial}{\partial t}f_{2,2}
&=& f_{2,2} + e^t\frac{y(1-y)}{x(1-x)}\cdot \frac{1}{2}x(1-x)\frac{\partial^2}{\partial x^2}f_{0,0}
\nonumber \\
&=& 
f_{2,2} + \frac{1}{2}\frac{\partial^2}{\partial x^2}x(1-x)f_{2,2}
\nonumber \\
&=& \frac{1}{2}x(1-x)\frac{\partial^2}{\partial x^2}f_{2,2}
+ \frac{1}{2}(2-4x)\frac{\partial}{\partial x}f_{2,2}.
\nonumber \\
%\label{calc:0}
\end{eqnarray*}
This shows that the generator in the $h$-transformed process is
\[
\frac{1}{2}x(1-x)\frac{\partial^2}{\partial x^2} +
\frac{1}{2}(2-4x)\frac{\partial}{\partial x},
\]
which corresponds to a Wright-Fisher diffusion with $\theta_1=\theta_2=2$.
The bridge density (\ref{bridge:rev:0}) is invariant under the $h$-transform $h(x)=x(1-x)$ and also under $h(x,t)=e^tx(1-x)$.
This choice of $h(x,t)$ is made because $x(1-x)$ is the first eigenfunction and the first eigenvalue is 1. 

There is a $h$-transform identity similar to (\ref{h2:0}) when $\theta_1=0,\theta_2=\theta >0$ where
\begin{equation}
				f_{2,\theta}(x,y;t) = e^{\theta t/2}f_{0,\theta}(x,y;t)y/x.
\label{h2:theta}
\end{equation}
A straightforward calculation is that
\begin{equation*}
\frac{\partial }{\partial t}f_{2,\theta} =
\frac{1}{2}x(1-x)\frac{\partial^2}{\partial x^2}f_{2,\theta}
+\frac{1}{2}\big (2 - (2+\theta)x\big )\frac{\partial}{\partial x}f_{2,\theta}.
%\label{h2:theta:1}
\end{equation*}
This shows that the generator in the $h$-transformed process is
\[
\frac{1}{2}x(1-x)\frac{\partial^2}{\partial x^2} +
\frac{1}{2}(2-(2+\theta)x)\frac{\partial}{\partial x},
\]
which corresponds to a Wright-Fisher diffusion with $\theta_1=2$, $\theta_2=\theta$.
The first eigenfunction in the process is $x$ and first eigenvalue $\theta/2$.
\begin{cor}
\label{cor:a}
An alternative form for the density of $X^{x,z,[0,T]}(t)$ in a neutral model when $\theta_1=\theta_2=0$, for $0 \leq x,z \leq 1$ is
\begin{eqnarray}
f_{x,z,[0,T]}(y;t) &=&\frac{f_{2,2}(x,y;t)f_{2,2}(y,z;T-t)}{f_{2,2}(x,z;T)}
\nonumber \\
&=& \frac{f_{2,2}(z,y;T-t)f_{2,2}(y,x;t)}{f_{2,2}(z,x;T)}
\nonumber \\
&=& f_{z,x,[0,T]}(y;T-t).
\label{bridge:rev:1}
\end{eqnarray}
An alternative form when $\theta_1=0,\theta_2=\theta>0$ $0 \leq x,z \leq 1$ is
\begin{eqnarray}
f_{x,z,[0,T]}(y;t) &=&\frac{f_{2,\theta}(x,y;t)f_{2,\theta}(y,z;T-t)}{f_{2,\theta}(x,z;T)}
\nonumber \\
&=& \frac{f_{2,\theta}(z,y;T-t)f_{2,\theta}(y,x;t)}{f_{2,\theta}(z,x;T)}
\nonumber \\
&=& f_{z,x,[0,T]}(y;T-t).
\label{bridge:rev:2}
\end{eqnarray}
\end{cor}
The densities on the right of these expressions are well defined at the boundaries when $x,z \in \{0,1\}$ because of non-zero recurrent mutation rates. 
Notice that the two alternative forms of the bridge density (\ref{bridge:rev:1}) and 
(\ref{bridge:rev:2}) are not continuous in $\theta$ at $\theta=0$.

The alternative forms become important in numerical calculations for the bridge density in a later section. 
\begin{cor} The limit density of $X^{x,z,[0,T]}(T/2+t)$, $\theta_1=0,\theta_2=0$ as $T \to \infty$ is independent of $t$ 
and is 
%the same as the quasi-stationary density of,
\begin{equation}
%\lim_{t\to \infty}
%\frac{P_x\big (X(t) \in (y,y+dy)\big )}{P_x\big (X(t) \in (0,1)\big )}
%=
 6y(1-y),\>0 < y < 1.
\label{quasibr:0}
\end{equation}
for $0 \leq x,z \leq 1$.
\end{cor}

\begin{cor} The limit density of $X^{x,z,[0,T]}(T/2+t)$, $\theta >0$ as $T \to \infty$ is independent of $t$ 
and is 
%the same as the quasi-stationary density of,
\begin{equation}
%\lim_{t\to \infty}
%\frac{P_x\big (X(t) \in (y,y+dy)\big )}{P_x\big (X(t) \in (0,1)\big )}
%=
 \theta y(1-y)^{\theta-1},\>0 < y < 1.
\label{quasibr:0}
\end{equation}
for $0 \leq x,z \leq 1$.
\end{cor}
The density (\ref{bridge:rev:1}) in Corollary \ref{cor:a} gives an alternative expression to (\ref{thm1:0}) by setting $x=z=0$, noting that terms on the right side are well defined in this case. 

Other identities for the coalescent transition functions follow.

Letting $x\to 0$ in the identity (\ref{h2:0}), using (\ref{fexp:3})  
\begin{eqnarray}
%&&
\sum_{r=0}^\infty q_r^4(t)(r+3)(r+2)y(1-y)^{r+1}
%\nonumber \\
%&&~~~~
= e^ty(1-y)\sum_{r=2}^\infty r(r-1)q_r^0(t)(1-y)^{r-2},
\end{eqnarray}
and equating coefficients of $y(1-y)^{l-1}$ shows the identity
\begin{equation}
(l+1)q_{l-2}^4(t) = e^t(l-1)q_l^0(t).
\label{h2:1}
\end{equation}
Another identity from the spectral expansion (\ref{fexp:1}), and (\ref{h2:0}), equating coefficients of $e^{-\frac{1}{2}n(n+4-1)t}$, is that for $n \geq 0$,
\begin{equation}
x(1-x)\widetilde{P}^{(1,1)}_n(r) = \widetilde{P}^{(-1,-1)}_{n+2}(r).
\label{Pidentity:0}
\end{equation}

Letting $x\to 0$ in the identity (\ref{h2:theta}) with $\theta_1=0,\theta_2=\theta>0$
\begin{eqnarray*}
&&
\sum_{r=0}^\infty q_r^{2+\theta}(t)(r+1+\theta)(r+\theta)y(1-y)^{r+\theta-1}
\nonumber \\
&&~~~~
= e^{\theta t/2}y\sum_{l=2}^\infty l(l+\theta -1)q_l^\theta(t)(1-y)^{l+\theta -2}
\end{eqnarray*}
and equating coefficients of $y(1-y)^{l+\theta-2}$ shows the identity
\begin{equation}
(l+\theta)q_{l-1}^{2+\theta}(t) = e^{\theta t/2}lq_l^\theta(t).
\label{h2:theta:4}
\end{equation}
Also from (\ref{fexp:1}) and (\ref{h2:theta}) 
\begin{equation}
x\widetilde{P}^{(\theta-1,1)}_n(r) = \frac{n+1}{n+\theta}\widetilde{P}_{n+1}^{(\theta-1,-1)}(r).
\label{Pidentity:2}
\end{equation}
There is interest in the way the $h$-transforms are used to show the pairs of identities (\ref{h2:1}), (\ref{h2:theta:4}) and (\ref{Pidentity:0}) and (\ref{Pidentity:2}). The first set is new, but probably not the second as so much is already known about Jacobi polynomials. The $h$-transform method of proof is new.
%%%
\subsection{Wright-Fisher bridges and branching P\'olya urns}\label{3.3}
%\textcolor{red}{Explain where this is going better.}\\
 A new representation of neutral Wright-Fisher bridges that come from Polya urns is now detailed. This has interest as a very classical probabilistic connection.

Genealogy in the Wright-Fisher diffusion process can be modelled by branching P{\'o}lya urns \citep{GS2010}. 
The joint density of $X=X(0)$ and $Y=X(t)$ when the marginal density of $X$ is Beta $(\theta_1,\theta_2)$, from the conditional density of $Y$ given $X$,  (\ref{fexp:3}), can be written as
\begin{eqnarray}
&&\sum_{l=0}^\infty q_l^\theta(t)\sum_{k=0}^l
{l\choose k}\frac{
{\theta_1}_{(k)}{\theta_2}_{(l-k)}
}
{
\theta_{(l)}
}
{\cal B}_{k+\theta_1,l-k+\theta_2}(x){\cal B}_{k+\theta_1,l-k+\theta_2}(y).
%\nonumber \\
%&&~~~~~~~~\times
%B(k+\theta_1,l-k+\theta_2)^{-1}x^{\theta_1+k-1}(1-x)^{\theta_2+l-k-1}
%\nonumber \\
%&&~~~~~~~~\times
%B(k+\theta_1,l-k+\theta_2)^{-1}y^{\theta_1+k-1}(1-y)^{\theta_2+l-k-1}.
\label{Polya:0}
\end{eqnarray}
There is a P{\'o}lya urn interpretation of (\ref{Polya:0}), explained in \citet{GS2010,GS2013}. 
Consider two P\'olya urns with an initial configuration of $\theta_1$ red and $\theta_2$ blue balls. The urns are identical until there are $l$ extra  balls, then branch to two urns with probability $q^\theta_l(t)$.
 The composition of both urns at the moment of splitting has a Beta binomial distribution
\[
{l\choose k}
\frac{
{\theta_1}_{(k)}
{\theta_2}_{(l-k)}
}{
\theta_{(l)}
},\>k=0,1,\ldots ,l.
\]
Draws continue independently in each urn after branching.  The unconditional relative  frequencies of red balls in the urns after an infinity of draws are Beta $(\theta_1,\theta_2)$ and the relative frequencies conditional on a configuration of $(k,l-k)$ are Beta  $(\theta_1+k,\theta_2+l-k)$. 
 In both cases the relative frequency distributions are the de Finetti measures for the sequences of draws.
The transition density $f_{\theta_1,\theta_2}(x,y;t)$ is the density of the relative frequency $Y$ in the second urn given the relative frequency $x$ in the first urn.
Instead of arguing to obtain a joint distribution of $X$ and $Y$ it is possible to argue directly to obtain the conditional distribution of $Y$ given $X=x$. 
If $X=x$, then since the distribution of $X$ is the de Finetti measure for whether the draws are red or blue, the probability of a configuration of $(k,l-k)$ balls given the urns branch after $l$ draws and $X=x$ is ${l\choose k}x^k(1-x)^{l-k}$.  the distribution of $Y$ given the $(k,l-k)$ configuration is then Beta $(\theta_1+l,\theta_2+l-k)$. Therefore the conditional distribution of $Y$ given $X=x$ is (\ref{fexp:3}).
 The conditional density $f_{\theta_1,\theta_2}(x,y;t)$ converges as $\theta_1,\theta_2 \to 0$ to
\begin{equation}
f_{0,0}(x,y;t) =	\sum_{l=1}^\infty q_l^0(t)\sum_{k=1}^{l-1}{l\choose k}x^k(1-x)^{l-k}{\cal B}_{k,l-k}(y).
%B(k,l-k)^{-1}y^{k-1}(1-y)^{l-k-1},\>0 < y < 1.
\label{density_zero:0}
\end{equation}
Relating the limit to convergence in an urn model: if $l>1$ the conditional distribution of the urn configuration at $l$ given $0 < x < 1$ as $\theta_1,\theta_2 \to 0$ is
\[
{l\choose k}x^k(1-x)^{l-k},\> 0 < k <l
\]
as a consequence of the de Finetti representation.  The joint probability that $k=0$ (and similarly for $k=l$) at trial $l$ and $0 < x < 1$ tends to zero as $\theta_1,\theta_2 \to 0$, even though it formally appears to be $(1-x)^l$.

We now give a new interpretation of Wright-Fisher bridges in terms of an urn model with three urns.
In a neutral Wright-Fisher diffusion with mutation the density of the frequency at time $t$ in a bridge from $x$ to $z$ at time $t$ is (\ref{bridgetr;0}), where $f_{\theta_1,\theta_2}(x,y;t)$ is the density
(\ref{fexp:0}) or (\ref{fexp:3}). An urn representation is found by considering the joint distribution of $X,Y,Z$ when $X$ and $Z$ have marginal Beta $(\theta_1,\theta_2)$ distributions. The joint distribution of the three random variables is then
\begin{eqnarray}
&&{\cal B}_{\theta_1,\theta_2}(x)f_{\theta_1,\theta_2}(x,y;t)f_{\theta_1,\theta_2}(y,z;T-t)
\nonumber \\
&=&
\sum_{l,m=0}^\infty q_l^\theta(t)q_m^\theta(T-t)
\sum_{j\leq l,k \leq m}
{n\choose j+k}
\frac{
{\theta_1}_{(j+k)}{\theta_2}_{(n-j-k)}
}
{
\theta_{(n)}
}
\frac {
{l\choose j}{m\choose k}
}
{
{n\choose j+k}
}
\nonumber \\
&&~~\times
{\cal B}_{\theta_1+j,\theta_2+l-j}(x){\cal B}_{\theta_1+k,\theta_2+m-k}(z)
{\cal B}_{\theta_1+k+j,\theta_2+n-k-j}(y),
\label{threeurn:0}
\end{eqnarray}
where $n=l+m$.
An urn representation has three urns $U,V,W$ which share draws before branching after draw $n$ in $W$.
The probability of branching after $n$ draws and that there are $l$, $m$ draws respectively in $U,V$ is $q_l(t)q_m(T-t)$.  The $l,m$ draws are chosen at random from the $n$ draws. After branching the urn draws continue independently in $U,V,W$. For fixed $l,m$ the probability of configurations $(j,l-j)$ and $(k,m-k)$ of (red,blue) balls in urns $U,V$ and $(j+k,n-j-k)$ in urn $W$ is
\[
{n\choose j+k}
\frac{
{\theta_1}_{(j+k)}{\theta_2}_{(n-j-k)}
}
{
\theta_{(n)}
}
\frac {
{l\choose j}{m\choose k}
}
{
{n\choose j+k}
}.
\]
Conditional on $l,m,j,k$ the limit relative frequencies of red balls in $U,V,W$ have the joint density in the last line of (\ref{threeurn:0}). The unconditional joint density of the relative frequencies in urns $U,V,W$ is the full expression (\ref{threeurn:0}).

%The conditional density of $X,Z$ given $Y=y$ as $\theta_1,\theta_2\to 0$  is
%\[
%f_{0,0}(y,x;t)f_{0,0}(y,z;T-t)
%\]
%because
% of reversibility and (\ref{density_zero:0}).
The configuration in the urns at branching, conditional on $X=x,Z=z$ has a probability distribution
\begin{eqnarray}
&&p_{\theta_1,\theta_2}(l,m;j,k;x,z;t) =
\nonumber \\
&&
\frac{
q_l^\theta(t)q_m^\theta(T-t)
{n\choose j+k}
\frac{
{\theta_1}_{(j+k)}{\theta_2}_{(n-j-k)}
}
{
\theta_{(n)}
}
\frac {
{l\choose j}{m\choose k}
}
{
{n\choose j+k}
}
{\cal B}_{\theta_1+j,\theta_2+l-j}(x){\cal B}_{\theta_1+k,\theta_2+m-k}(z)
}
{
{\cal B}_{\theta_1,\theta_2}(x)f_{\theta_1,\theta_2}(x,z;T)
}.
\nonumber \\
\label{cond_urn:1}
\end{eqnarray}
Then for $l,m\geq 1$ as $\theta_1,\theta_2 \to 0$ the distribution (\ref{cond_urn:1}) converges to 
\begin{eqnarray}
&&p_{0,0}(l,m;j,k;x,z;t) =
\nonumber \\
&&
\frac{
q_l^0(t)q_m^0(T-t)
{n\choose j+k}
\frac{
{(j+k-1)!}{(n-j-k-1)!}
}
{
(n-1)!
}
\frac {
{l\choose j}{m\choose k}
}
{
{n\choose j+k}
}
{\cal B}_{j,l-j}(x){\cal B}_{k,m-k}(z)
}
{
x^{-1}(1-x)^{-1}f_{0,0}(x,z;T)
}.
\nonumber \\
\label{cond_urn:2}
\end{eqnarray}
There is a similar limit if $\theta_1 \to 0$ and $\theta_2 > 0$.
The relative frequency of red balls in urn $W$ after an infinite number of draws when $\theta_1,\theta_2 >0$ has a Beta $(\theta_1+j+k,\theta_2+n-j-k)$ distribution so the density of $Y$ conditional on $X=x,Z=z$ is
\begin{equation}
\sum_{l,m=0}^\infty p_{\theta_1,\theta_2}(l,m;j,k;x,z;t){\cal B}_{\theta_1+j+k,\theta_2+n-j-k}(y),
\label{cond_urn:3}
\end{equation}
with a similar expression when $\theta_1=0,\theta_2=0$ where summation is over $l,m \geq 1$, $1 \leq  j  \leq  l-1$ and $1 \leq k \leq m-1$. The density (\ref{cond_urn:3}) is the same as the neutral bridge densities in Section \ref{NWFB} with the different cases of mutation. 
The limit distribution when $\theta_1=0,\theta_2=0$ and $x,z \to 0$ is now considered. The denominator of (\ref{cond_urn:2}) converges to 
\[
\sum_{l=1}^\infty q_l^0(T)l(l-1)
\]
which simplifies as (\ref{ffact:0}).  The conditional probability distribution of the urn configurations 
(\ref{cond_urn:2}) converges to zero unless $j=1,k=1$ when the limit is
\begin{eqnarray*}
%&&
p_{0,0}(l,m;j,k;0,0;t) &=&
%\nonumber \\
%&&
\frac{
q_l^0(t)q_m^0(T-t)
{n\choose 2}
\frac{
{1!}{(n-3)!}
}
{
(n-1)!
}
\frac {
lm
}
{
{n\choose 2}
}
(l-1)(m-1)
}
{
\sum_{l=2}^\infty e^{-\frac{1}{2}l(l-1)T}(2l-1)(l-1)
}
\nonumber \\
&=&
\frac{
q_l^0(t)q_m^0(T-t)
\frac {
l(l-1)m(m-1)
}
{
(n-1)(n-2)
}
}
{
\sum_{l=2}^\infty e^{-\frac{1}{2}l(l-1)T}(2l-1)(l-1)
}.
%\label{cond_urn:4}
\end{eqnarray*}
The conditional distribution of the relative frequency of red balls in urn $W$ is then
\begin{equation*}
\sum_{n=2}^\infty \sum_{l,m\geq 1,l+m=n}p_{0,0}(l,m;j,k;0,0;t){\cal B}_{2,n-2}(y).
%\label{cond_urn:5}
\end{equation*}
%%%
%%%
%%
\subsection{Wright-Fisher bridges with genic selection}\label{3.4}
In this section we consider the genealogy of a Wright-Fisher bridge when there is selection in the model. This is a new approach different from that in \citet{SGE2013}. The genealogy of the Wright-Fisher diffusion with selection is more complex that in a neutral model and the transition functions for the coalescent genealogy do not have an explict form.

Theorem \ref{thm:a} holds when $\gamma\ne 0$, if $x \ne 0$, however care is again needed when beginning from $x=0$. For definiteness let $\gamma >0$,
then the speed measure density when $\theta_1=0$ and $\theta_2 = \theta \geq 0$ is
\[
e^{\gamma x}x^{-1}(1-x)^{\theta-1},\>0 < x < 1.
\]
The density is reversible before fixation and with $\bm{\theta}=(0,\theta)$,
\begin{eqnarray}
f_{0,\theta,\gamma}(x,y;t) 
&=& 
\sum_{\bm{\alpha}} b^{(0,\theta)}_{\bm{\alpha}}(t,x)\pi[\bm{\alpha}+\bm{\theta}](y)
\nonumber \\
&=& 
\sum_{\bm{\alpha}} b^{(0,\theta)}_{\bm{\alpha}}(t,y)\pi[\bm{\alpha}+\bm{\theta}](x)
\frac{y^{-1}(1-y)^{\theta-1}e^{\gamma y}}
{x^{-1}(1-x)^{\theta-1}e^{\gamma x}}.
\label{rev:gamma}
\end{eqnarray}
$\alpha_i \geq 1$ in the summation, $i=1,2$. In a neutral model
\[
b^{(0,\theta)}_{\bm{\alpha}}(t,y) = {|\bm{\alpha}|\choose \alpha_1}y^{\alpha_1}(1-y)^{\alpha_2}
\]
An asymptotic form as $x \to 0$ is
\begin{equation*}
f_{0,\theta,\gamma}(x,y;t)=y^{-1}(1-y)^{\theta-1}e^{\gamma y}x 
\sum_{\alpha_2\geq 1} b^{(0,\theta)}_{(1,\alpha_2+\theta)}(t,y)(\alpha_2+\theta)
c(1,\alpha_2+\theta)^{-1} + o(x^2),
%\label{asymptotic:0}
\end{equation*}
since $\pi[\bm{\alpha}+\bm{\theta}](x) = o(x^{\alpha_1-1})$ as $ x \to 0$ for $\alpha_1 > 1$ and 
\[
\pi[(1,\alpha_2+\theta)](x) = c(1,\alpha_2+\theta)^{-1}\frac{\Gamma(1+\alpha_2+\theta)}{\Gamma(1)\Gamma(\alpha_2+\theta)}x^{1-1}(1-x)^{\alpha_2+\theta-1}
e^{\gamma x}.
\] 
Note that
\begin{equation*}
c(1,\alpha_2+\theta) =
\int_0^1(\alpha_2+\theta)z^{\alpha_2+\theta-1}e^{z\gamma}dz,\>\gamma \geq  0.
\end{equation*}
Denote
\[
g_\theta(y;t;\gamma) = 
y^{-1}(1-y)^{\theta-1}e^{\gamma y}
\sum_{\alpha_2=1}^\infty b^{(0,\theta)}_{(1,\alpha_2)}(t,y)(\alpha_2+\theta)c(1,\alpha_2+\theta)^{-1}.
\]
We have shown that as $x \to 0$
\[
f_{0,\theta,\gamma}(x,y;t)=xg_\theta(y;t;\gamma) + o(x^2)
\]
and that
\[
f_{0,\theta,\gamma}(y,0;t)=
y(1-y)^{1-\theta}e^{-\gamma y}
g_\theta(y;t;\gamma).
\]
The next theorem is an analogue of Theorem \ref{thm:b} which follows by substitution  in the second line of the bridge density (\ref{bridge:rev:0}) and the asymptotic  forms in this section for the transition density.

\begin{thm}\label{thm:c:gamma}
The density of $X^{0,0,[0,T]}(t)$, $0\leq t \leq T$ in a Wright-Fisher diffusion bridge when $\gamma >0$ with mutation away from $A$ at rate $\theta/2 \geq 0$ is
\begin{equation*}
\frac{
y(1-y)^{1-\theta}e^{-\gamma y}g_\theta(y;t,\gamma)g_\theta(y;T-t,\gamma)
}
{
\sum_{\alpha_2=1}^\infty {b^\prime}^{(0,\theta)}_{(1,\alpha_2)}(T,0)c(1,\alpha_2+\theta)^{-1}(\alpha_2+\theta)
},
%\label{thm:d:gamma}
\end{equation*}
where 
%\[
%g_\theta(y;t;\gamma) = \sum_{\alpha_2=1}^\infty b^{(0,\theta)}_{(1,\alpha_2)}(t,0)c(1,\alpha_2+\theta)^{-1}(\alpha_2+\theta)
%\]
%and
prime denotes differentiation with respect to $y$.
\end{thm}
%Theorems \ref{thm:b:gamma} and
The theorem shows that in a non-neutral $(0,0)$ Wright-Fisher bridge there is exactly one lineage from $t$ to $0$ and from $t$ to $T$ and that there is a reversibility from either end of the bridge. This is a simlar interpretation to in the neutral bridge, Theorems \ref{thm:b} and \ref{thm:c}. An  illustration of the genealogy is in Figure 4. Theorem \ref{thm:c:gamma} and the coalescent interpretation are new.

\begin{center}
\vbox{
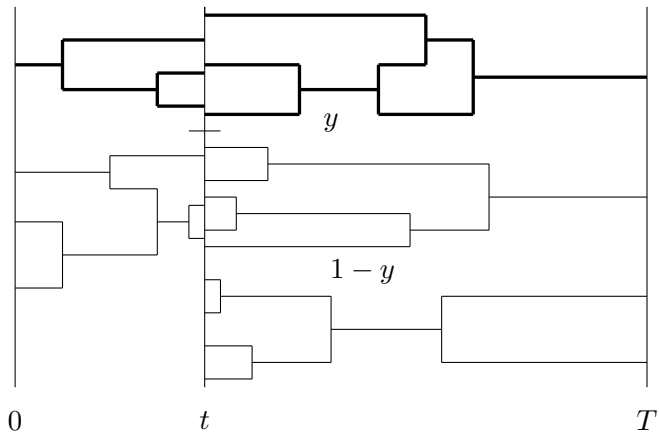
 Figure 4: Genealogy in a non-neutral (0,0) bridge
\bigskip

\begin{tikzpicture}[xscale = 0.21,yscale=0.22]
%\draw (0,0) -- (0,21);
\draw (0,0) -- (0,23);
%\draw (12,0) -- (12,21);line width=1.25pt
\draw (12,0) -- (12,15.5);
%\draw [line width=1.25pt] (12,15.5) -- (12,21);
%\draw [line width=1.25pt] (12,15.5) -- (12,22.5);
\draw (12,15.5) -- (12,23);
\draw (12,22.5) -- (12,23);
%\draw (40,0) -- (40,21);
\draw (40,0) -- (40,23);
%%%Left side
\draw [line width=1.25pt] (0,19.5) -- (3,19.5);
\draw [line width=1.25pt] (3,21) -- (12,21);
\draw [line width=1.25pt] (3,18) -- (9,18);
\draw [line width=1.25pt] (9,19) -- (12,19);
\draw [line width=1.25pt] (9,17) -- (12,17);
\draw [line width=1.25pt] (3,21) -- (3,18);
\draw [line width=1.25pt] (9,19) -- (9,17);
\draw (0,13) -- (6,13);
\draw (6,14) -- (12,14);
\draw (6,12) -- (9,12);
\draw (0,10) -- (3,10);
\draw (0,6) -- (3,6);
\draw (3,8) -- (9,8);
\draw (9,10) -- (11,10);
\draw (11,11) -- (12,11);
\draw (11,9) --  (12,9);
\draw (3,6) -- (3,10);
\draw (6,14) -- (6,12);
\draw (9,12) -- (9,8);
\draw (11,11) -- (11,9);
%
%%
%\draw (0,13) -- (6,13);
%\draw (6,14) -- (12,14);
%\draw (6,12) -- (12,12);
%\draw (6,14) -- (6,12);
%%%
%\draw (0,8) -- (5,8);
%\draw (5,7) -- (12,7);
%\draw (5,10) -- (7,10);
%\draw (7,11) -- (12,11);
%\draw (7,9) -- (12,9);
%\draw (5,10) -- (5,7);
%\draw (7,11) -- (7,9);
%% Right side
%\draw [line width=1.25pt] (40,18) -- (18,18);
\draw [line width=1.25pt] (18,19.5) -- (12,19.5);
\draw [line width=1.25pt] (18,16.5) -- (12,16.5);
\draw [line width=1.25pt] (18,19.5) -- (18,16.5);
\draw [line width=1.25pt] (12,22.5) -- (26,22.5);
\draw [line width=1.25pt] (23,18) -- (18,18);
\draw [line width=1.25pt] (23,19.5) -- (23,16.5);
\draw [line width=1.25pt] (29,21) -- (29,16.5);
\draw [line width=1.25pt] (29,18.75) -- (40,18.75);
\draw [line width=1.25pt] (23,16.5) -- (29,16.5);
\draw [line width=1.25pt] (23,19.5) -- (26,19.5);
\draw [line width=1.25pt] (26,21.0) -- (29,21.0);
\draw [line width=1.25pt] (26,22.5) -- (26,19.5);
\draw (40,11.5) -- (30,11.5);
\draw (30,13.5) -- (16,13.5);
\draw (16,14.5) -- (12,14.5);
\draw (16,12.5) -- (12,12.5);
\draw (16,14.5)  -- (16,12.5);
\draw (30,13.5) -- (30,9.5);
\draw (25,8.5) -- (12,8.5);
\draw (25,8.5) -- (25,10.5);
\draw (14,11.5) -- (14,9.5);
\draw (30,9.5) -- (25,9.5);
\draw (25,10.5) -- (14,10.5);
\draw (14,9.5) -- (12,9.5);
\draw (14,11.5) -- (12,11.5);
\draw (27,3.5) -- (20,3.5);
\draw (27,1.5) -- (27,5.5);
\draw (27,1.5) -- (40,1.5);
\draw (27,5.5) -- (40,5.5);

\draw (20,5.5) -- (13,5.5);
\draw (13,4.5) -- (12,4.5);
\draw (13,6.5) -- (12,6.5);
\draw (13,6.5) -- (13,4.5);
\draw (20,1.5) -- (15,1.5);
\draw (15,0.5) -- (12,0.5);
\draw (15,2.5) -- (12,2.5);
\draw (15,2.5) -- (15,0.5);
\draw (20,5.5) -- (20,1.5);
%%
%\draw (-1,14.5) -- (1,14.5);
%\draw -- (-3,17.25) node {$x$};
%\draw -- (-3,7.25) node {$1-x$};
%%
%\draw (39,16) -- (41,16);
%\draw -- (43,18) node {$z$};
%\draw -- (43,8.75) node {$1-z$};
%%
\draw (11,15.5) -- (13,15.5);
\draw -- (20,16) node {$y$};
\draw -- (22,7) node {$1-y$};
\draw -- (0,-2) node {$0$};
\draw -- (12,-2) node {$t$};
\draw -- (40,-2) node {$T$};
\end{tikzpicture}
}
\end{center}
%\bigskip

We are interested in finding the Yaglom density and the limit distribution of $X^{x,z,[0,T]}(T/2+t)$ as $T\to \infty$. These calculations depend on the first eigenfunction $w(y)$. The Yaglom density is proportional to $y^{-1}(1-y)^{\theta-1}e^{\gamma y}w(y)$.
An analogy of the $h$-transform (\ref{h2:0}) in the non-neutral model is 
\[
f_w(x,y;t) = f(x,y;t)\frac{w(y)}{w(x)}e^{\lambda t},
\]
where $(\lambda,w)$ are the first eigenvalue, eigenfunction pair satisfying (\ref{eigenfunction:0}) when either or both of $\theta_1,\theta_2$ are zero. A calculation similar to that in the neutral model shows that the generator associated with the transition functions $f_w$ when $\theta_1=0,\theta_2=\theta \geq 0$ is
\begin{equation}
				\frac{1}{2}x(1-x)\frac{\partial^2}{\partial x^2} + \frac{1}{2}\Bigg (\gamma x(1-x) + 2x(1-x)\frac{\partial}{\partial x}\log w(x) - \theta x\Bigg )\frac{\partial}{\partial x}.
\label{fw:0}
\end{equation}
This is a generator of a Wright-Fisher diffusion with frequency dependent selection. The stationary distribution of the process with generator (\ref{fw:0}) is  
\[
\psi (y) = Cw(y)^2\pi(y) = Cw(y)^2e^{\gamma y}y^{-1}(1-y)^{\theta-1}.
\]
An alternative form for the density of $X^{x,z,[0,T]}$ is, for $0 \leq t \leq T$
\begin{equation}
\frac{f_w(x,y;t)f_w(y,z;T-t)}{f_w(x,z;T)}.
\label{wbridge:0}
\end{equation}
The limit density as $T \to \infty$ of $X^{x,z,[0,T]}(T/2+t)$, from 
(\ref{wbridge:0}) is 
\[
\psi(y)\psi(z)/\psi(z) = \psi(y).
\]
Computing $\psi(y)$ (when $\theta=0$) rests on finding a numerical solution of the eigenfunction equation
\begin{equation*}
y(1-y)w^{\prime\prime}(y) + \gamma y(1-y)w^\prime(y) + 2\lambda w(y) = 0
%\label{wde:0}
\end{equation*}
for a maximal $\lambda$. There is a unique solution with $\lambda >0, w(x) \geq 0$. Clearly $w(0)=w(1)=0$ but the maximal $\lambda$ is unknown.
Parameterizing $\lambda \equiv \lambda(\gamma)$, then $\lambda (0) = 1$, the first eigenvalue in a neutral model. \cite{K1955} recognizes a connection with Spheroidal wave functions, tabulates $\lambda (\gamma)$, and expresses $w(x)$ as a series in Gegenbauer polynomials. The spheroidal wave equations are 
\begin{equation}
\frac{d}{dz}\left ((1-z^2)\frac{dS}{dz}\right ) +
\left (\mu + c^2(1-z^2) - \frac{m^2}{1-z^2}\right )S = 0,
\label{Sph:0}
\end{equation}
see \citet{NIST2010}.
If $c^2 > 0$ the equations are \emph{prolate} or $c^2 < 0$ \emph{oblate}.
Solutions to (\ref{Sph:0}) are the Spheroidal wave functions $P_{S_n^m}(z,c^2)$ with eigenvalues $\mu_n^m(c^2)$, $n=m,m+1,\ldots$. 
Letting 
\[
w(y) = e^{-(\gamma/2)y}y^{1/2}(1-y)^{1/2}S^\circ(y),
\] 
then
\begin{equation*}
\frac{d}{dy}\left (y(1-y)\frac{dS^\circ}{dy}\right )
+ \left (\mu - c^2 + 4c^2y(1-y) - \frac{m^2}{4y(1-y)}\right )S^\circ = 0.
%\label{Sph:1}
\end{equation*}
This shows that $e^{\gamma/2y}y^{-1/2}(1-y)^{-1/2}w(y)$ is proportional to 
$P_{S_1^1}(2y-1),-\gamma^2/16)$. Parameters are $m=1$, $c^2 = -\gamma^2/16$, and $\mu - c^2 = 2\lambda$ and $z=1-2y$. Scaling $w(y)$ to be a probability distribution 
\begin{equation*}
w(y) =
 \frac{
 e^{-\gamma/2y}
 y^{1/2}(1-y)^{1/2}
 P_{S_1^1}(2y-1,-\gamma^2/16)
 }
 {
 \int_0^1e^{-(\gamma/2)x}
 x^{1/2}(1-x)^{1/2}
 P_{S_1^1}(2x-1,-\gamma^2/16)dx
 }.
%\label{wsoln:0}
\end{equation*}
Spheroidal wave functions are implemented in Mathematica and it is straightforward to calculate and plot $w(y)$ via WolframAlpha on a web interface. The basic command for solutions to (\ref{Sph:0}) is {\tt SpheroidalPS[n,m,c,z]}. In our application $c^2 \leq 0$ so $c$ is a multiple of the complex number $i$.
For example if $\gamma=3$, then $c=0.75i$, and $n=m=1$, 
\begin{verbatim}
Plot[(-0.166599)^(-1)*exp(-1.5*y)*(y*(1-y))^(1/2)
 *SpheroidalPS[1,1,0.75i,2*y-1], {y,0,1}]
\end{verbatim}
is the appropriate command. The constant $-0.166599$ is the scaling factor returned by
\begin{verbatim}
Integrate[exp(-1.5*x)(x*(1-x))^(1/2)
 *SpheroidalPS[1,1,0.75i,2*x -1],{x,0,1}]
\end{verbatim}
It is not necessary to find the first eigenvalue for the above commands, however {\tt SpheroidalEigenvalue[n,m,c]} will return it. In this example
{\tt Eigenvalue[1,1,0.75i]} returns $2.44853$, so $\lambda=1.505515$.
The full expansion of the transition function is 
\begin{eqnarray*}
&&f_{0,0,\gamma}(x,y;t) \nonumber = e^{\gamma y}y^{-1}(1-y)^{-1}
\sum_{n=1}^\infty e^{-(\lambda^1_n(\frac{\gamma^2}{16})+\frac{\gamma^2}{16}) \frac{t}{2}}\frac{2n+1}{n(n+1)}
\nonumber \\
&&~~~~~~~~~
\times e^{-\frac{\gamma}{2}x}x^{\frac{1}{2}}(1-x)^{\frac{1}{2}}
P_{S^1_n}(2x-1,-\frac{\gamma^2}{16})
\nonumber \\
&&~~~~~~~~~ \times e^{-\frac{\gamma}{2}y}y^{\frac{1}{2}}(1-y)^{\frac{1}{2}}
P_{S^1_n}(2y-1,-\frac{\gamma^2}{16}).
%\nonumber \\
%\label{swvexp:0}
\end{eqnarray*}
This expansion is analogous to the neutral expansion (\ref{fexp:2}). 
Mathematica also includes a function for the derivative of a Spheroidal wave function, {\tt SpheroidalPSPrime[n,m,c,z]}. The drift coefficient in the $h$-transformed generator (\ref{fw:0}) when $\theta=0$ can thus be easily evaluated or plotted. The coefficient is easily seen to be
\[
\frac{1}{2}(1-2x) + 2x(1-x)
\frac{P^\prime_{S_1^1}(2x-1,-\gamma^2/16)}
{P_{S_1^1}(2x-1,-\gamma^2/16)}.
\]
Carrying through the example above, if $\gamma=3$, the command to plot the drift coefficient is
\begin{verbatim}
Plot[0.5-x +2.0*x*(1-x)*SpheroidalPSPrime[1,1,0.75i,2*x-1]
 /SpheroidalPS[1,1,0.75i,2*x-1],{x,0,1}]
\end{verbatim}
In this particular example the the plot shows that the drift coefficient is very close to $1-2x$, which would be exact if $\gamma=0$ and $w(x) = x(1-x)$.
\subsection*{Appendix: Simulating the neutral Wright-Fisher bridge}
It is possible to simulate exactly from the density of a Wright-Fisher bridge in a neutral model. We develop a new algorithm here when one or both of $x,z$ are zero. \citet{JS2015} gave a sophisticated algorithm for this in the special case that both $\theta_1,\theta_2 > 0$ and $x,z > 0$. The identities in this paper
% \comment{PJ}{i.e.\ Bob's identities!} 
provide a new direct method for other cases considered here:
\begin{enumerate}
\item $\theta_1=0$, $\theta_2 = \theta > 0$, and $x = z = 0$.
\item $\theta_1=\theta_2 = 0$ and $x = z = 0$.
\item $\theta_1=0$, $\theta_2 = \theta > 0$, and $x, z > 0$.
\item $\theta_1=\theta_2 = 0$ and $x, z > 0$.
\end{enumerate}
\subsection*{1. $\theta_1=0$, $\theta_2 = \theta > 0$ and $x = z = 0$}
The density of $X^{0,0,[0,T]}(t)$ is given in Theorem \ref{thm:c}.
%\footnote{Terms in red refer to Bob's draft of 11/1/2016.}.
 It can be written as
\begin{equation}
\label{eq:mixture}
f^{0,0,[0,T]}(t) = \sum_{l_1 = 1}^\infty \sum_{l_2 = 1}^\infty p_{l_1,l_2} \frac{y(1-y)^{l_1+l_2+\theta-3}}{B(2,l_1+l_2+\theta-2)},
\end{equation}
where
%\footnote{$h^\theta(T)$ is the denominator of \red{(45)}.}
\begin{align*}
p_{l_1,l_2} &= \frac{1}{h^\theta(T)}\frac{l_1(l_1+\theta-1)l_2(l_2+\theta-1)}{(l_1+l_2+\theta-1)(l_1+l_2+\theta-2)}q_{l_1}^\theta(t)q_{l_2}^\theta(T-t),\\
h^\theta(T) &= g_\theta(0;T) = \sum_{l=1}^\infty e^{-l(l+\theta-1)T/2}(2l+\theta - 1)l(l+\theta-1).
\end{align*}
We recognise \eqref{eq:mixture} as the density of an infinite mixture of Beta random variables, so that $(p_{l_1,l_2})_{l_1,l_2 \in \bbN}$ is a probability mass function on $\bbN^2$ (for convenience, set $p_{l_1,l_2} = 0$ if $l_1=0$ or $l_2 =0$). The $(l_1,l_2)$th component corresponds to a $\text{Beta}(2,l_1+l_2+\theta-2)$ variate. Therefore in order to simulate from $f^{0,0,[0,T]}(t)$, it suffices to simulate from the discrete distribution $(p_{l_1,l_2})$. This is complicated by the infinite series representations for $q_{l_1}^\theta(t)$, $q_{l_2}^\theta(T-t)$, and $h^\theta(T)$, but a solution is possible via the \emph{series method} \citep[Ch.\ IV.5]{D1986}; we summarise the strategy as follows.

Suppose $(p_l)_{l\in\bbN}$ is a probability mass function whose masses may not be computable with finite resource but for which we have available a pair of sequences $(p^-_l(k))_{k\in\bbN}$, $(p^+_l(k))_{k\in\bbN}$ for each $l$ with the following properties:
\begin{enumerate}%[(i)]
\item $p^-_l(k) \leq p_l \leq p^+_l(k)$ for all $l,k \in \mathbb{N}$.
\item For each $l$, $p^-_l(k) \uparrow p_l$ as $k \to\infty$.
\item For each $l$, $p^+_l(k) \downarrow p_l$ as $k \to\infty$.
\end{enumerate}
For $U \sim \text{Uniform}[0,1]$, standard inversion sampling implies that $\inf\{ L \in \bbN: \sum_{l=0}^L p_l > U\}$ is distributed according to $(p_l)$, but this is not computable, as we have noted. However, it is easily verified that if
\[
K(L) = \inf\left\{k\in \bbN: \sum_{l=0}^L p^-_l(k) > U \text{ or } \sum_{l=0}^L p^+_l(k) < U\right\}
\]
then
\[
\inf\left\{ L \in \bbN: \sum_{l=0}^L p^-_l(K(l)) > U \right\}
\]
is also distributed as $(p_l)$, and this can be computed from finitely many terms in the double arrays $(p^-_l(k))$, $(p^+_l(k))$.

To employ this strategy for $(p_{l_1,l_2})$ in \eqref{eq:mixture} we must find monotonically converging upper and lower bounds on each of $q_{l_1}^\theta(t)$, $q_{l_2}^\theta(T-t)$, and $h^\theta(T)$. The first two expressions are covered by \citet[Proposition 1]{JS2015}, who showed that
\begin{multline*}
(q^\theta)^-_l(k,t) := \sum_{j=l}^{l+2k+1} \rho_j^\theta(t)(-1)^{j-l}\frac{(2j+\theta-1)(l+\theta)_{j-1}}{l!(j-l)!} \leq q^\theta_l(t) \\
\leq  \sum_{j=l}^{l+2k} \rho_j^\theta(t)(-1)^{j-l}\frac{(2j+\theta-1)(l+\theta)_{j-1}}{l!(j-l)!} =: (q^\theta)^+_l(k,t),
\end{multline*}
provided $2k+l \geq C^{(t,\theta)}_l := \inf\left\{i \geq l: \frac{\theta + i+l-1}{i-l+1}\frac{\theta + 2i+1}{\theta + 2i-1}e^{-(i + \theta/2)t} < 1\right\}$, a constant beyond which convergence of these bounds becomes monotonic in $k$. \citet{JS2015} assumed that $\gq_1, \gq_2 > 0$, but their result follows without change when either or both of $\gq_1 = 0$, $\gq_2 = 0$.] It remains to find similar converging bounds on $h^\theta(T)$. Writing
\[
h^\theta(T) = \sum_{l=1}^\infty e^{-l(l+\theta-1)T/2}(2l+\theta - 1)l(l+\theta-1) =: \sum_{l=1}^\infty h^\theta_l(T),
\]
we find for $l\geq 2$ that
\[
h^\theta_{l+1}(T) = \frac{(2l+1+\theta)(l+1)(l+\theta)}{(2l-1+\theta)(l-1)(l+\theta-1)}e^{-(l+\theta/2)T}h^\theta_l(T) \leq 5e^{-(l+\theta/2)T}h^\theta_l(T) < h^\theta_l(T),
\]
with the final inequality holding provided $l > D^{(\theta,T)} := \frac{\log 5}{T} - \frac{\theta}{2}$. It follows immediately that, for $k > \max(2,D^{(\theta,T)})$,
\[
0 < \sum_{l=k}^\infty h^\theta_l(T) \leq \sum_{l=k}^\infty (5e^{-(k+\theta/2)T})^{l-k}h^\theta_k(T) = \frac{h^\theta_k(T)}{1-5e^{-(k+\theta/2)T}},
\]
and hence that
\begin{eqnarray*}
				h^-_k(T) &:=& \sum_{l=1}^k e^{-l(l+\theta-1)T/2}(2l+\theta - 1)l(l+\theta-1),\nonumber \\
				h^+_k(T) &:=& h^-_k(T) + \frac{h^\theta_{k+1}(T)}{1-5e^{-(k+1 + \theta/2)T}}
\end{eqnarray*}
serve as monotonically converging lower and upper bounds on $h^\theta(T)$ as $k\to\infty$. We summarise the above argument in the following result.

\begin{thm}
Let $\Sigma : \bbN \to \bbN^2$ be any bijective pairing function denoted by $\Sigma(l)=(l_1,l_2)$, and let
\begin{align*}
S^-_\bfk(L) &= \ds\sum_{l=\Sigma(0)}^{\Sigma(L)}\frac{l_1(l_1+\theta-1)l_2(l_2+\theta-1)}{(l_1 + l_2 +\theta- 1)(l_1 + l_2 +\theta - 2)} \frac{(q^\theta)_{l_1}^-(k_l,t)(q^\theta)_{l_2}^-(k_l,T-t)}{h^+_{k_l}(T)},\\
S^+_\bfk(L) &= \ds\sum_{l=\Sigma(0)}^{\Sigma(L)}\frac{l_1(l_1+\theta-1)l_2(l_2+\theta-1)}{(l_1 + l_2 +\theta- 1)(l_1 + l_2 +\theta- 2)} \frac{(q^\theta)_{l_1}^+(k_l,t)(q^\theta)_{l_2}^+(k_l,T-t)}{h^-_{k_l}(T)}.
\end{align*}
If $k_l > \max\left((C_{l_1}^{(t)}-l_1)/2, (C_{l_2}^{(T-t)}-l_2)/2, D^{(T,\theta)},2\right)$ for each $l=0,1,\dots, L$, then $S^-_\bfk(L)$ and $S^+_\bfk(L)$ are respectively lower and upper bounds on the distribution function $\sum_{l=\Sigma(0)}^{\Sigma(L)} p_{l_1,l_2}$ which are monotonically converging to this limit as $\bfk \to\infty$ (that is, as $k_l \to \infty$ for every component of $\bfk = (k_0,k_1,\dots,k_L)$). Thus the following algorithm returns exact samples from the density of $X^{0,0,[0,T]}(t)$.

\begin{algorithm}[H]
\DontPrintSemicolon
Set $l \longleftarrow 0$, $k_0 \longleftarrow 0$, $\bfk \longleftarrow (k_0)$.\;
Simulate $U \sim \text{Uniform}[0,1]$.\;
\Repeat{false}{
Set $(l_1,l_2) \longleftarrow \Sigma(l)$, $k_l \longleftarrow \left\lceil \max\left((C_{l_1}^{(t)}-l_1)/2, (C_{l_2}^{(T-t)}-l_2)/2, D^{(T,\theta)},2\right)\right\rceil$.\;
\While{$S_\bfk^-(l) < U < S_\bfk^+(l)$}{Set $\bfk \longleftarrow \bfk + (1,1,\ldots,1)$.\nllabel{q6}}\;
\uIf{$S_\bfk^-(l) > U$}{\Return{$Y \sim \text{Beta}(2,l_1+l_2+\theta-2)$.}}
\ElseIf{$S_\bfk^+(l) < U$}{Set $\bfk \longleftarrow (k_0,k_1,\ldots,k_l, 0)$.\;
Set $l \longleftarrow l+1$.}
}
\caption{Simulating from the density $f^{0,0,[0,T]}(t)$ of $X^{0,0,[0,T]}(t)$ when $\theta_1 = 0$, $\theta_2 = \theta > 0$.}%\label{alg:EAWF0q00}
\end{algorithm}
\end{thm}

\subsection*{2. $\theta_1=\theta_2 = 0$ and $x = z = 0$}
This case is very similar to that of the previous subsection and so is omitted. (In fact, one can go through the previous subsection and replace $\theta$ with $0$ throughout.) The density of the bridge when $\theta_1 = \theta_2 = 0$ simplifies slightly (see Theorem \ref{thm:b}) to
\begin{eqnarray*}
%%\label{eq:mixture2}
%\ds
 f^{0,0,[0,T]}(t) &=& \sum_{l_1 = 2}^\infty \sum_{l_2 = 2}^\infty p_{l_1,l_2} \frac{y(1-y)^{l_1+l_2-3}}{B(2,l_1+l_2-2)},
\nonumber \\
 %\ds 
p_{l_1,l_2} &=& \frac{1}{h(T)}\frac{l_1(l_1-1)l_2(l_2-1)}{(l_1+l_2-1)(l_1+l_2-2)}q_{l_1}^0(t)q_{l_2}^0(T-t).
\end{eqnarray*}

\subsubsection*{3. $\theta_1=0$, $\theta_2 = \theta > 0$ and $x, z > 0$}
By Corollary \ref{cor:a}, the density of $X^{x,z,[0,T]}(t)$ in this case is the same as for $\theta_1 = 2$, $\theta_2 = \theta$, which is covered by Algorithm 4 of \citet{JS2015}.
\subsection*{4. $\theta_1=\theta_2 = 0$ and $x, z > 0$}
By Corollary \ref{cor:a}, the density of $X^{x,z,[0,T]}(t)$ in this case is the same as for $\theta_1 = \theta_2 = 2$, which is covered by Algorithm 4 of \citet{JS2015}.

\end{document}